\documentclass[a4paper,11pt,reqno]{amsart}

\usepackage{amssymb}

\input prepictex
\input pictex
\input postpictex

\def\SunoScuatrodos{\beginpicture
 \setcoordinatesystem units <1pt, 1pt>
 \multiput{$\circ$} at 0 0  60 0  /
 \put{$\bullet$} at 30 0
 \put{{\tiny $\times$}} at 60 0
 \put{{\small$\alpha_1$}} at 0 6
 \put{{\small $\alpha_2$}} at 30 6
 \put{{\small$\alpha_3$}} at 60 6
 \putrule from 2 1 to 28 1
 \putrule from 2 -1 to 28 -1
 \putrule from 32 1 to 58 1
 \putrule from 32 -1 to 58 -1
 \put{$>$} at 15 0
 \put{$<$} at 45 0
 \endpicture}

\def\ScuatroScuatrodos{\beginpicture
 \setcoordinatesystem units <1pt, 1pt>
 \multiput{$\circ$} at 0 0 20 0 40 0 60 0  75 10 75 -10 /
 \multiput{{\tiny$\times$}} at 60 0 75 -10 /
 \putrule from 2 0 to 18 0
 \putrule from 22 0 to 38 0
 \putrule from 42 0 to 58 0
 \plot  62 1  73 9 /
 \plot  62 -1  73 -9 /
 \put{{\small$\alpha_1$}} at 0 6
 \put{{\small$\alpha_2$}} at 20 6
 \put{{\small$\alpha_3$}} at 40 6
 \put{{\small$\alpha_4$}} at 58 6
 \put{{\small$\alpha_5$}} at 83 10
 \put{{\small$\alpha_6$}} at 83 -10
 \endpicture}

\def\SochoScuatrodos{\beginpicture
 \setcoordinatesystem units <1pt, 1pt>
 \multiput{$\circ$} at 0 0 20 0 40 0 60 0 80 0 100 0 40 -20 /
 \multiput{{\tiny$\times$}} at 80 0 100 0 /
 \putrule from 2 0 to 18 0
 \putrule from 22 0 to 38 0
 \putrule from 42 0 to 58 0
 \putrule from 62 0 to 78 0
 \putrule from 82 0 to 98 0
 \putrule from  40 -2 to 40 -18
 \put{{\small$\alpha_1$}} at 0 6
 \put{{\small$\alpha_2$}} at 48 -20
 \put{{\small$\alpha_3$}} at 20 6
 \put{{\small$\alpha_4$}} at 40 6
 \put{{\small$\alpha_5$}} at 60 6
 \put{{\small$\alpha_6$}} at 80 6
 \put{{\small$\alpha_7$}} at 100 6
 \endpicture}

\def\ScuatrodosScuatrodos{\beginpicture
 \setcoordinatesystem units <1pt, 1pt>
 \multiput{$\circ$} at 0 0 20 0 40 0 60 0  75 10 75 -10 /
 \multiput{{\tiny$\times$}} at 0 0 20 0 60 0 75 -10 /
 \putrule from 2 0 to 18 0
 \putrule from 22 0 to 38 0
 \putrule from 42 0 to 58 0
 \plot  62 1  73 9 /
 \plot  62 -1  73 -9 /
 \setdots<2pt>
 \put{{\small$\alpha_1$}} at 0 6
 \put{{\small$\alpha_2$}} at 20 6
 \put{{\small$\alpha_3$}} at 40 6
 \put{{\small$\alpha_4$}} at 58 6
 \put{{\small$\alpha_5$}} at 83 10
 \put{{\small$\alpha_6$}} at 83 -10
 \endpicture}

\def\SunodosScuatrodos{\beginpicture
 \setcoordinatesystem units <1pt, 1pt>
 \multiput{$\circ$} at 0 0 20 0 60 0  /
 \put{$\bullet$} at 40 0
 \multiput{{\tiny $\times$}} at 0 0 60 0 /
 \putrule from 2 0 to 18 0
 \putrule from 22 1 to 38 1
 \putrule from 22 -1 to 38 -1
 \putrule from 42 1 to 58 1
 \putrule from 42 -1 to 58 -1
 \put{$>$} at 30 0
 \put{$<$} at 50 0
 \put{{\small$\alpha_1$}} at 0 6
 \put{{\small $\alpha_2$}} at 20 6
 \put{{\small$\alpha_3$}} at 40 6
 \put{{\small$\alpha_4$}} at 60 6
 \endpicture}

\def\SdosScuatrodos{\beginpicture
 \setcoordinatesystem units <1pt, 1pt>
 \multiput{$\circ$} at 0 0 20 0 40 0  55 10 55 -10 /
 \multiput{{\tiny$\times$}} at 40 0  55 -10 /
 \putrule from 2 0 to 18 0
 \putrule from 22 0 to 38 0
 \plot  42 1  53 9 /
 \plot  42 -1  53 -9 /
 \put{{\small$\alpha_1$}} at 0 6
 \put{{\small$\alpha_2$}} at 20 6
 \put{{\small$\alpha_3$}} at 38 6
 \put{{\small$\alpha_4$}} at 63 10
 \put{{\small$\alpha_5$}} at 63 -10
 \endpicture}

\def\SunodosSunodos{\beginpicture
 \setcoordinatesystem units <1pt, 1pt>
 \multiput{$\circ$} at 0 0 20 0 40 0 60 0  /
 \put{{\tiny $\times$}} at  60 0
 \putrule from 2 0 to 18 0
 \putrule from 22 1 to 38 1
 \putrule from 22 -1 to 38 -1
 \putrule from 42 0 to 58 0
 \put{$>$} at 30 0
 \put{{\small$\alpha_1$}} at 0 6
 \put{{\small $\alpha_2$}} at 20 6
 \put{{\small$\alpha_3$}} at 40 6
 \put{{\small$\alpha_4$}} at 60 6
 \endpicture}

\def\SunoSunodos{\beginpicture
 \setcoordinatesystem units <1pt, 1pt>
 \multiput{$\circ$} at 0 0 30 0 60 0  /
 \put{{\tiny $\times$}} at 30 0
 \putrule from 2 0 to 28 0
 \putrule from 32 0 to 58 0
 \put{{\small$\alpha_1$}} at 0 6
 \put{{\small $\alpha_2$}} at 30 6
 \put{{\small$\alpha_3$}} at 60 6
 \endpicture}

\def\ScuatroSunodos{\beginpicture
 \setcoordinatesystem units <1pt, 1pt>
 \multiput{$\circ$} at 0 0 30 0 60 0 30 -25 /
 \put{{\tiny $\times$}} at 30 -25
 \putrule from 2 0 to 28 0
 \putrule from 32 1 to 58 1
 \putrule from 32 -1 to 58 -1
 \putrule from 30 -2 to 30 -23
 \put{$<$} at 45 0
 \put{{\small$\alpha_1$}} at 0 6
 \put{{\small $\alpha_2$}} at 30 6
 \put{{\small$\alpha_3$}} at 60 6
  \put{{\small$\alpha_4$}} at 38 -25
 \endpicture}

\def\SochoSunodos{\beginpicture
 \setcoordinatesystem units <1pt, 1pt>
 \multiput{$\circ$} at 0 0 20 0 40 0 60 0 80 0   /
 \put{{\tiny$\times$}} at 80 0
 \putrule from 2 0 to 18 0
 \putrule from 22 1 to 38 1
 \putrule from 22 -1 to 38 -1
 \putrule from 42 0 to 58 0
 \putrule from 62 0 to 78 0
 \put{$>$} at 30 0
 \put{{\small$\alpha_1$}} at 0 6
 \put{{\small$\alpha_2$}} at 20 6
 \put{{\small$\alpha_3$}} at 40 6
 \put{{\small$\alpha_4$}} at 60 6
 \put{{\small$\alpha_5$}} at 80 6
 \endpicture}

\def\SdosSunodos{\beginpicture
 \setcoordinatesystem units <1pt, 1pt>
 \multiput{$\circ$} at 0 0 40 0 20 20    /
 \put{{\tiny$\times$}} at 20 20
 \putrule from 2 0 to 38 0
 \plot 2 1 18 19 /
 \plot 38 1 22 19 /
 \put{{\small$\alpha_1$}} at 0 -6
 \put{{\small$\alpha_2$}} at 40 -6
 \put{{\small$\alpha_3$}} at 28 20
 \endpicture}

\newtheorem{theorem}[equation]{Theorem}
\newtheorem{lemma}[equation]{Lemma}
\newtheorem{corollary}[equation]{Corollary}
\newtheorem{proposition}[equation]{Proposition}

\numberwithin{equation}{section}

\theoremstyle{definition}

\newtheorem*{example*}{Example}
\newtheorem{example}[equation]{Example}

\newtheorem{remark}[equation]{Remark}
\newtheorem*{remark*}{Remark}

\newcommand{\bN}{{\mathbb N}}
\newcommand{\bZ}{{\mathbb Z}}

\newcommand{\frg}{{\mathfrak g}}

\newcommand{\frh}{{\mathfrak h}}
\newcommand{\fri}{{\mathfrak i}}
\newcommand{\frc}{{\mathfrak c}}
\newcommand{\frd}{{\mathfrak d}}
\newcommand{\frf}{{\mathfrak f}}
\newcommand{\fre}{{\mathfrak e}}
\newcommand{\frt}{{\mathfrak t}}

\newcommand{\frb}{{\mathfrak b}}

\newcommand{\frn}{{\mathfrak n}}

\newcommand{\frgoo}{{\frg_{(\bar 0,\bar 0)}}}
\newcommand{\frgouno}{{\frg_{(\bar 0,\bar 1)}}}
\newcommand{\frgunoo}{{\frg_{(\bar 1,\bar 0)}}}
\newcommand{\frgunouno}{{\frg_{(\bar 1,\bar 1)}}}
\newcommand{\calS}{{\mathcal S}}

\newcommand{\subo}{_{\bar 0}}
\newcommand{\subuno}{_{\bar 1}}

\providecommand{\espan}[1]{\text{span}\left\{ #1\right\}}

 \DeclareMathOperator{\tri}{\mathfrak{tri}}

 \DeclareMathOperator{\frosp}{\mathfrak{osp}}
 \DeclareMathOperator{\frsl}{{\mathfrak{sl}}}
  
 \DeclareMathOperator{\frsp}{{\mathfrak{sp}}}
 \DeclareMathOperator{\frso}{{\mathfrak{so}}}
 \DeclareMathOperator{\frpsl}{{\mathfrak{psl}}}
 \DeclareMathOperator{\frgl}{{\mathfrak{gl}}}
 \DeclareMathOperator{\frpgl}{{\mathfrak{pgl}}}

 \DeclareMathOperator{\frj}{{\mathfrak{j}}}

 \DeclareMathOperator{\End}{End}
 
 \DeclareMathOperator{\Mat}{Mat}

 \DeclareMathOperator{\charac}{char}
  \DeclareMathOperator{\rank}{rank}

\def\bigstrut{\vrule height 12pt width 0ptdepth 2pt}
\def\hregla{\hrule height.1pt}
\def\hreglabis{\hrule height .3pt depth -.2pt}
\def\hregleta{\hrule height .5pt}
\def\hreglon{\hrule height1pt}
\def\vreglon{\vrule height 12pt width1pt depth 4pt}
\def\vregleta{\vrule width .5pt}
\def\hreglonfill{\leaders\hreglon\hfill}
\def\hreglafill{\leaders\hreglabis\hfill}
\def\hregletafill{\leaders\hregleta\hfill}
\def\vregla{\vrule width.1pt}

\newenvironment{romanenumerate}
 {\begin{enumerate}
 
 }{\end{enumerate}}

\begin{document}

\title{An extended Freudenthal Magic Square in characteristic $3$}

\author[Isabel Cunha]{Isabel Cunha$^{\diamond}$}
 \thanks{$^{\diamond}$ Supported by Centro de Matem\'atica da Universidade de Coimbra -- FCT}
 \address{Departamento de Matem\'atica, Universidade da Beira
 Interior,\newline 6200 Covilh\'a, Portugal}
 \email{icunha@mat.ubi.pt}

\author[Alberto Elduque]{Alberto Elduque$^{\star}$}
 \thanks{$^{\star}$ Supported by the Spanish Ministerio de
 Educaci\'on y Ciencia
 and FEDER (MTM 2004-081159-C04-02) and by the
Diputaci\'on General de Arag\'on (Grupo de Investigaci\'on de
\'Algebra)}
 \address{Departamento de Matem\'aticas, Universidad de
Zaragoza, 50009 Zaragoza, Spain}
 \email{elduque@unizar.es}

\date{May 15, 2006}



\begin{abstract}
Freudenthal's Magic Square, which in characteristic $0$ contains the
exceptional Lie algebras other than $G_2$, is extended over fields
of characteristic $3$, through the use of symmetric composition
superalgebras, to a larger square that contains both Lie algebras
and superalgebras. With one exception, the simple Lie superalgebras
that appear have no counterpart in characteristic $0$.
\end{abstract}

\maketitle


\section{Introduction}

In 1966 \cite{Tits66}, Tits gave a unified construction of the
exceptional simple Lie algebras based on two ingredients: a unital
composition algebra (or Hurwitz algebra) and a central simple degree
$3$ Jordan algebra. This construction is valid over arbitrary fields
of characteristic $\ne 2,3$. The outcome of this construction is
given by Freudenthal's Magic Square \cite{Schaferbook,Freudenthal64}
(Table \ref{ta:Freudenthalsquare}).

\begin{table}[h!]
$$ \vbox{\offinterlineskip
 \halign{$#$&\hfil\quad$#$\quad\hfil&%
 \vreglon #
 &\hfil\quad$#$\quad\hfil&\hfil\quad$#$\quad\hfil
 &\hfil\quad$#$\quad\hfil&\hfil\quad$#$\ \hfil\cr
 &&\omit&\multispan4{\hfil$\dim J$\hfil}\cr
 \noalign{\smallskip}
 \bigstrut &&width 0pt&6&9&15&27\cr
 &&\multispan5{\hreglonfill}\cr
 &1&&A_1&A_2&C_3&F_4\cr
 \bigstrut&2&&A_2&\omit$A_2\oplus A_2$&A_5&E_6\cr
 \bigstrut\smash{\raise 6pt\hbox{$\dim A$}}&4&&C_3&A_5&D_6&E_7\cr
 \bigstrut&8&&F_4&E_6&E_7&E_8\cr}}
$$
\bigskip
\caption{Freudenthal's Magic Square}\label{ta:Freudenthalsquare}
\end{table}

At least in the split cases, this is a construction that depends on
two Hurwitz algebras, since the central simple degree Jordan
algebras turn out to be the algebras of hermitian $3\times 3$
matrices over a Hurwitz algebra. Even though the construction is not
symmetric on the two Hurwitz algebras involved, the result (the
Magic Square) is symmetric.

\smallskip

Over the years, several symmetric constructions of  Freudenthal's
Magic Square based on two Hurwitz algebras have been proposed.
Vinberg, in 1966 \cite{Vinberg05} gave a construction based on two
Hurwitz algebras and their derivation algebras. In 1993, Allison and
Faulkner \cite{AllisonFaulkner93} gave a very general construction
based on structurable algebras. In particular, the tensor product of
two Hurwitz algebras is structurable and, in this case, the
construction is equivalent to the one given by Vinberg. More
recently, Barton and Sudbery \cite{BS00,BS03} and Landsberg and
Manivel \cite{LM02,LM04} provided a different construction based on
two Hurwitz algebras and their Lie algebras of triality. (See
\cite{Baez} for a very nice exposition.)

The triality phenomenon was described by E.~Cartan \cite{Cartan38},
relating the natural and the two half spin representations of the
orthogonal Lie algebra in dimension $8$. As shown in \cite{KMRT},
simpler formulas can be obtained if the so called \emph{symmetric
composition algebras} are used, instead of the classical Hurwitz
algebras. This led the second author to reinterpret the construction
by Barton and Sudbery in terms of a construction of a Lie algebra
$\frg(S,S')$ out of two symmetric composition algebras $S$ and $S'$
\cite{ElduqueMagicI}. With a few exceptions in dimension $2$, any
symmetric composition algebra of dimension $1$, $2$ or $4$ is a
\emph{para-Hurwtiz} algebra, while in dimension $8$, besides the
para-Hurwitz algebras, there appears a new family of symmetric
composition algebras, the \emph{Okubo} algebras. The construction in
\cite{ElduqueMagicI} using para-Hurwitz algebras reduces naturally
to the construction by Barton and Sudbery (although with slightly
simpler formulas). Okubo algebras provide new constructions that
highlight different order $3$ automorphisms and different
distinguished subalgebras of the excepctional simple Lie algebras.
Further descriptions of the Lie algebras that appear in
Freudenthal's Magic Square in terms of copies of the
three-dimensional simple Lie algebra $\frsl_2$ and its two
dimensional natural module are derived in \cite{ElduqueMagicII}.

\smallskip

An interesting advantage of the constructions of the Magic Square in
terms of two composition algebras and their triality Lie algebras is
that this construction is also valid over fields of characteristic
$3$. Only some attention has to be paid to the second row (or
column). Thus \cite{ElduqueMagicI}, for instance, instead of a
simple Lie algebra of type $E_6$, a non simple Lie algebra of
dimension $78$ is obtained if a two dimensional and an eight
dimensional composition algebras are used as ingredients, but it
contains a simple codimension 1 ideal of type $E_6$ (the simple Lie
algebra of type $E_6$ has dimension $77$ in characteristic $3$).

The characteristic $3$ presents also another exceptional feature.
Only over fields of this characteristic there are composition
superalgebras which are nontrivial (that is, they are not
composition algebras). In \cite{EldOkuSuperCompo} both Hurwitz
superalgebras and symmetric composition superalgebras are
classified. The main results assert that the only nontrivial
superalgebras appear in characteristic $3$ and dimensions $3$ and
$6$, in agreement with the classification of the simple alternative
superalgebras by Shestakov \cite{Shestakov97}.

The symmetric composition superalgebras can thus be plugged into the
construction $\frg(S,S')$ in \cite{ElduqueMagicI}. Then an extension
in characteristic $3$ is obtained of  Freudenthal's Magic Square, in
which Lie superalgebras appear (see Table
\ref{ta:superFreudenthal}).

\begin{table}[h!]
$$
\vbox{\offinterlineskip
 \halign{$#$\qquad&\hfil$#$\quad\hfil&%
 \vreglon #%
 &\hfil\quad$#$\quad\hfil&\hfil\quad$#$\quad\hfil
 &\hfil\quad$#$\quad\hfil&\hfil\quad$#$\quad\hfil&%
 \vrule  depth 4pt width .5pt #%
 &\hfil\quad$#$\quad\hfil&\hfil\quad$#$\ \hfil\cr
 &&\omit&\multispan7{\hfil$\dim S$\qquad\hfil}\cr
 \noalign{\medskip}
 \bigstrut &&width 0pt&1&2&4&8&\omit\vrule height 8pt depth 4pt width .5pt&3&6\cr
 &&\multispan8{\hreglonfill}\cr
 &1&&A_1&\tilde A_2&C_3&F_4&&(6,8)&(21,14)\cr
 \bigstrut&2&& &\omit$\tilde A_2\oplus \tilde A_2$&\tilde A_5&\tilde E_6&&(11,14)&(35,20)\cr
 \bigstrut\smash{\lower 6pt\hbox{$\dim S'$}}&4&& & &D_6&E_7&&(24,26)&(66,32)\cr
 \bigstrut&8&& & & &E_8&&(55,50)&(133,56)\cr
 &\multispan9{\hregletafill}\cr
 \bigstrut& 3&& & & & & & (21,16)&(36,40)\cr
 \bigstrut& 6&& & & & & & & (78,64)\cr}}
$$
\bigskip
\caption{Freudenthal's Magic Supersquare (characteristic
$3$)}\label{ta:superFreudenthal}
\end{table}

In Table \ref{ta:superFreudenthal} only the entries above the
diagonal have been displayed, as the square is symmetric. The
notation $\tilde X$ indicates that in characteristic $3$ the algebra
obtained is not simple, but contains a simple codimension $1$ ideal
of type $X$. For instance, in the split case, $\tilde A_5$ denotes
the projective general Lie algebra $\frpgl_6$. Besides, only the
pairs $(\dim\frg\subo,\dim\frg\subuno)$ are displayed for the
nontrivial Lie superalgebras that appear in the ``supersquare''.

\smallskip

The aim of this paper is the description of all the Lie
superalgebras that appear in Table \ref{ta:superFreudenthal}. With
just a single exception, which corresponds to the pair $(6,8)$,
these Lie superalgebras have no counterpart in characteristic $0$,
and they are either simple or contain a simple codimension $1$ ideal
(this happens again only in the second row). They will be described
as contragredient Lie superalgebras, and to do so, a previous
description in terms similar to those used in \cite{ElduqueMagicII}
will be provided too.

In a forthcoming paper, most of these Lie superalgebras will be
shown to be isomorphic to the Lie superalgebras constructed in terms
of orthogonal and symplectic triple systems in
\cite{ElduqueNewSimple3}.

\smallskip

The paper is organized as follows. Section 2 will be devoted to
review the symmetric composition superalgebras and their Lie
superalgebras of triality. Then Section 3 will deal with the
construction of a Lie superalgebra in terms of two symmetric
composition superalgebras, thus extending the construction in
\cite{ElduqueMagicI}. Section 4 will give the definitions and
characterizations of contragredient Lie superalgebras in a way
suitable for our purposes. Finally, the last Section 5 will be
devoted to the description of all the split Lie superalgebras in
Freudenthal's Magic Supersquare as contragredient Lie superalgebras.
Their even and odd parts will be computed too.

\bigskip

\section{Symmetric composition superalgebras}

Recall \cite{EldOkuSuperCompo} that a \emph{quadratic superform} on
a $\bZ_2$-graded vector space $U=U\subo\oplus U\subuno$ over a field
$k$ is a pair $q=(q\subo,b)$ where
\begin{romanenumerate}
\item $q\subo :U\subo\rightarrow k$ is a quadratic form.

\item $b:U\times U\rightarrow k$ is a supersymmetric even bilinear form.
That is, $b\vert_{U\subo\times U\subo}$ is symmetric,
$b\vert_{U\subuno\times U\subuno}$ is alternating
($b(x\subuno,x\subuno)=0$ for any $x\subuno\in U\subuno$) and
$b(U\subo,U\subuno)=0=b(U\subuno,U\subo)$.

\item $b\vert_{U\subo\times U\subo}$ is the polar of $q\subo$. That
is,
$b(x\subo,y\subo)=q\subo(x\subo+y\subo)-q\subo(x\subo)-q\subo(y\subo)$
for any $x\subo,y\subo\in U\subo$.
\end{romanenumerate}

The quadratic superform $q=(q\subo,b)$ is said to be \emph{regular}
if $q\subo$ is regular (definition as in \cite[p.~xix]{KMRT}) and
the restriction of $b$ to $U\subuno$ is nondegenerate.

\smallskip

Then a superalgebra $C=C\subo\oplus C\subuno$ over $k$, endowed with
a regular quadratic superform $q=(q\subo,b)$, called the
\emph{norm}, is said to be a \emph{composition superalgebra} in case
\begin{subequations}\label{eq:norm}
\begin{align}
&q\subo(x\subo,y\subo)=q\subo(x\subo)q\subo(y\subo),\\
&b(x\subo y,x\subo z)=q\subo(x\subo)b(y,z)=b(yx\subo,zx\subo),\\
&b(xy,zt)+(-1)^{\lvert x\rvert\lvert y\rvert +\lvert x\rvert\lvert
 z\rvert+\lvert y\rvert\lvert z\rvert}b(zy,xt)=(-1)^{\lvert
 y\rvert\lvert z\rvert}b(x,z)b(y,t),
\end{align}
\end{subequations}
for any $x\subo,y\subo\in C\subo$ and homogeneous elements
$x,y,z,t\in C$ (where $\lvert x\rvert$ denotes the parity of the
homogeneous element $x$).

The unital composition superalgebras are termed \emph{Hurwitz
superalgebras}.

\smallskip

A composition superalgebra $S=S\subo\oplus S\subuno$ with norm
$q=(q\subo,b)$ is said to be \emph{symmetric} if
\[
b(xy,z)=b(x,yz),
\]
for any $x,y,z\in S$. That is, the associated supersymmetric
bilinear form is associative.

\smallskip

\begin{example}\label{ex:B12}
Let $V$ be a two dimensional vector space over the field $k$ with a
nonzero alternating bilinear form $\langle .\vert .\rangle$.
Consider the superspace
\[
B(1,2)=k1\oplus V,
\]
where $B(1,2)\subo=k1$ and $B(1,2)\subuno =V$, with supercommutative
multiplication given by
\[
1x=x1=x\qquad\text{and}\qquad uv=\langle u\vert v\rangle 1
\]
for any $x\in B(1,2)$ and $u,v\in V$, and with quadratic superform
$q=(q\subo,b)$ given by:
\[
q\subo(1)=1,\quad b(u,v)=\langle u\vert v\rangle,
\]
for any $u,v\in V$.

If the characteristic of $k$ is $3$, then $B(1,2)$ is a Hurwitz
superalgebra \cite[Proposition 2.7]{EldOkuSuperCompo}.
\end{example}

\medskip

\begin{example}\label{ex:B42}
Let $V$ be as in Example \ref{ex:B12}, then $\End_k(V)$ is equipped
with the symplectic involution $f\mapsto \bar f$, such that
\[
\langle f(u)\vert v\rangle =\langle u\vert\bar f(v)\rangle,
\]
for any $u,v\in V$. Consider the superspace
\[
B(4,2)=\End_k(V)\oplus V,
\]
where $B(4,2)\subo=\End_k(V)$ and $B(4,2)\subuno=V$, with
multiplication given by:
\begin{itemize}
\item the usual multiplication (composition of maps) in
$\End_k(V)$,
\smallskip

\item $v\cdot f=f(v)=\bar f\cdot v$ for any $f\in \End_k(V)$ and
$v\in V$,
\smallskip

\item $u\cdot v=\langle .\vert u\rangle v\,\bigl(w\mapsto
\langle w\vert u\rangle v\bigr)\,\in \End_k(V)$ for any $u,v\in V$,
\end{itemize}
\smallskip

\noindent and with quadratic superform $q=(q\subo,b)$ such that
\[
q\subo(f)=\det f,\quad b(u,v)=\langle u\vert v\rangle,
\]
for any $f\in \End_k(V)$ and $u,v\in V$.

Again, if the characteristic is $3$, $B(4,2)$ is a Hurwitz
superalgebra \cite[Proposition 2.7]{EldOkuSuperCompo}.

\medskip

The vector space $\End_k(V)$ may be identified with $V\otimes V$
(unadorned tensor products are assumed to be tensor products over
the ground field $k$) by means of
\[
\begin{split}
V\otimes V&\longrightarrow \End_k(V)\\
x\otimes y&\mapsto \langle x\vert .\rangle y:v\mapsto \langle x\vert
v\rangle y.
\end{split}
\]
The symplectic involution on $\End_k(V)$ becomes now
$\overline{x\otimes y}=-y\otimes x$ for any $x,y\in V$.

Then $B(4,2)$ can be identified with
\begin{equation}\label{eq:B42Vs}
B(4,2)=(V\otimes V)\oplus V,
\end{equation}
with multiplication given by
\[
\left\{\begin{aligned} &(x\otimes y)\cdot(z\otimes t)=\langle x\vert
t\rangle z\otimes y,\\
&u\cdot (x\otimes y)=\langle x\vert u\rangle y=-(y\otimes x)\cdot
u\\
&u\cdot v=-u\otimes v,
\end{aligned}\right.
\]
for any $x,y,z,t,u,v\in V$, and the supersymmetric bilinear form $b$
is determined by
\begin{equation}\label{eq:B42b}
b(x\otimes y,z\otimes t)=\langle x\vert z\rangle\langle y\vert
t\rangle,\quad b(u,v)=\langle u\vert v\rangle,
\end{equation}
for any $x,y,z,t,u,v\in V$ (see \cite[\S 2]{ElduqueMagicII}).
\end{example}

\medskip

Given any Hurwitz superalgebra $C$ with norm $q=(q\subo,b)$, the
linear map $x\mapsto \bar x=b(x,1)1-x$ is a superinvolution
($\bar{\bar x}=x$ and $\overline{xy}=(-1)^{\lvert x\rvert\lvert
y\rvert}\bar y\bar x$ for any homogeneous $x,y\in C$). Then $C$,
with the same norm $q$, but with new multiplication
\[
x\bullet y=\bar x\bar y
\]
becomes a symmetric composition superalgebra, which is called a
\emph{para-Hurwitz superalgebra}, and denoted by $\bar C$.

\smallskip

Also, as for the ungraded case, if $\varphi$ is an order $3$
automorphism of a symmetric composition superalgebra $S$, with
multiplication $\bullet$ and norm $q$ (we will refer to this as the
symmetric composition superalgebra $(S,\bullet,q)$), then $\varphi$
is an isometry of $q$, and with the new multiplication given by
\begin{equation}\label{eq:Sphi}
x*y=\varphi(x)\bullet\varphi^2(y),
\end{equation}
$(S,*,q)$ is again a symmetric composition algebra, denoted by
$S_\varphi$ or $(S,\bullet,q)_\varphi$.

\smallskip

\begin{example}\label{ex:B12bar}
Given any order $3$ automorphism of $B(1,2)$ over a field $k$ of
characteristic $3$ (and hence this is also an automorphism of the
associated para-Hurwitz superalgebra $\overline{B(1,2)}$), there is
a scalar $\lambda\in k$ and a symplectic basis $\{v,w\}$ of $V$
(that is, $\langle v\vert w\rangle =1$) such that $\varphi(v)=v$ and
$\varphi(w)=\lambda v +w$. The symmetric composition superalgebra
defined on $B(1,2)$ by means of the new multiplication given by
\[
x\bullet y=\varphi(\bar x)\varphi^2(\bar y)
\]
is denoted by $B(1,2)_\lambda$ (see \cite[Example
2.9]{EldOkuSuperCompo}).

For $\lambda =1$ this is just the para-Hurwitz superalgebra
$\overline{B(1,2)}$.
\end{example}

\medskip

\begin{example}\label{ex:B42bar}
In terms of the description of $B(4,2)$ as $(V\otimes V)\oplus V$,
the multiplication in the associated para-Hurwtiz superalgebra
$\overline{B(4,2)}$ is determined by:
\[
\left\{\begin{aligned} &(x\otimes y)\bullet(z\otimes t)=\langle
y\vert
z\rangle t\otimes x,\\
&u\bullet (x\otimes y)=\langle y\vert u\rangle x=-(y\otimes
x)\bullet
u\\
&u\bullet v=-u\otimes v,
\end{aligned}\right.
\]
for any $x,y,z,t,u,v\in V$. The associated bilinear form is given in
\eqref{eq:B42b}.
\end{example}

\medskip

\begin{remark} This is exactly the multiplication in the para-Cayley
algebra in \cite[eq. (2.2)]{ElduqueMagicII}, if the last copy of $V$
is suppressed there.
\end{remark}

\smallskip

The classification of the symmetric composition superalgebras
appears in \cite[Theorem 4.3]{EldOkuSuperCompo}. Over fields of
characteristic $2$, any symmetric composition superalgebra is a
symmetric composition algebra suitably graded over $\bZ_2$. In
characteristic $\ne 2$, the classification is given by:

\begin{theorem} Let $k$ be a field of characteristic $\ne 2$, and
let $S$ be a symmetric composition superalgebras over $k$. Then
either:
\begin{romanenumerate}
\item $S\subuno=0$; that is, $S$ is a symmetric composition algebra,
or
\item the characteristic of $k$ is $3$ and either there is a scalar
$\lambda\in k$ such that $S$ is isomorphic to
$\overline{B(1,2)}_\lambda$ or $S$ is isomorphic to
$\overline{B(4,2)}$.
\end{romanenumerate}
\end{theorem}

\smallskip

Therefore, the superalgebras in Examples \ref{ex:B12bar} and
$\ref{ex:B42bar}$ exhaust, up to isomorphism, all the symmetric
composition superalebras, which do not appear in the ungraded
setting.

\smallskip

In order to superize the construction of the Lie algebras
$\frg(S,S')$ in the Introduction, the triality Lie superalgebras are
needed. Given a symmetric composition superalgebra $(S,*,q)$, its
\emph{triality Lie superalgebra} $\tri(S,*,q)=\tri(S,*,q)\subo\oplus
\tri(S,*,q)\subuno$ is defined by
\[
\begin{split}
\tri&(S,*,q)_i\\
 &=\bigl\{ (d_0,d_1,d_2)\in \frosp(S,q)_i^3 :
d_0(x*y)=d_1(x)*y+(-1)^{i\lvert x\rvert}x*d_2(y)\\
 &\qquad\qquad\text{for any homogeneous elements $x,y\in S$}\bigr\},
\end{split}
\]
where $i=\bar 0,\bar 1$. For simplicity, this Lie superalgebra will
be referred to simply by $\tri(S)$. Here $\frosp(S,q)$ denotes the
orthosymplectic Lie superalgebra of the superform $q$.
The bracket in $\tri(S)$ is given
componentwise.

Note that $\tri(S)$ is endowed with a natural automorphism
$\theta:(d_0,d_1,d_2)\mapsto (d_2,d_0,d_1)$, satisfying
$\theta^3=1$.

\smallskip

\begin{theorem}\label{th:triality1224}
Let $(S,*,q)$ be any of the symmetric composition superalgebras
$\overline{B(1,2)}$ or $\overline{B(4,2)}$. Then for any homogeneous
element $d_0\in\frosp(S,q)_i$ ($i=\bar 0,\bar 1$), there are unique
$d_1,d_2\in\frosp(S,q)_i$ such that $(d_0,d_1,d_2)\in\tri(S)$.
Besides, the map $\Phi:\frosp(S,q)\rightarrow \frosp(S,q)$ given by
$\Phi(d_0)=d_1$ is an automorphism of $\frosp(S,q)$ with $\Phi^3=1$.
Hence $\tri(S)=\{\bigl(d,\Phi(d),\Phi^2(d)\bigr): d\in\frosp(S,q)\}$
is isomorphic to $\frosp(S,q)$. Moreover, for $\overline{B(1,2)}$,
this automorphism $\Phi$ is the identity.
\end{theorem}
\begin{proof} This is \cite[Theorem 5.6]{EldOkuSuperCompo}.
The last assertion is proved in \cite[p. 5464]{EldOkuSuperCompo}.
\end{proof}

\smallskip

\begin{corollary}\label{co:triB12bar}
$\tri(\overline{B(1,2)})=\{(d,d,d) : d\in \frosp(B(1,2),q)\}$.
\end{corollary}

\medskip

Given a nondegenerate even supersymmetric bilinear form $b$ on a
superspace $W=W\subo\oplus W\subuno$, its orthosymplectic Lie
superalgebra $\frosp(W,b)$ is spanned by the operators
\[
\sigma_{x,y}: z\mapsto (-1)^{\lvert y\rvert\lvert
z\rvert}b(x,z)y-(-1)^{\lvert x\rvert(\lvert y\rvert+\lvert
z\rvert)}b(y,z)x,
\]
for homogeneous $x,y,z\in W$. Besides, for homogeneous $x,y,z,t\in
W$,
\begin{equation}\label{eq:sigmas}
[\sigma_{x,y},\sigma_{z,t}]=\sigma_{\sigma_{x,y}(z),t} +
(-1)^{(\lvert x\rvert +\lvert y\rvert)\lvert
z\rvert}\sigma_{z,\sigma_{x,y}(t)}.
\end{equation}

Also, given a vector space $U$ endowed with a nondegenerate
alternating bilinear form $\langle .\vert .\rangle$, its symplectic
Lie algebra $\frsp\bigl(U,\langle .\vert.\rangle\bigr)$ is spanned
by the operators
\begin{equation}\label{eq:gammas}
\gamma_{x,y}:z\mapsto \langle x\vert z\rangle y+\langle y\vert
z\rangle x,
\end{equation}
for $x,y,z\in U$.

\medskip

Over any field $k$ of characteristic $3$, consider the para-Hurwitz
superalgebra $\overline{B(1,2)}$ which, for simplicity, will be
denoted by $S_{1,2}$. Thus, $S_{1,2}=k1\oplus V$, where $V$ is a two
dimensional vector space endowed with a nonzero (hence
nondegenerate) alternating bilinear form $\langle .\vert .\rangle$.
The norm $q=(q\subo,b)$ satisfies $q\subo(1)=1$ (so $b(1,1)=2=-1$)
and $b(u,v)=\langle u\vert v\rangle$ for any $u,v\in V$. The
corresponding orthosymplectic Lie superalgebra
$\frosp(S_{1,2},b)=\sigma_{S_{1,2},S_{1,2}}$ satisfies:
\[
\frosp(S_{1,2},b)\subo=\sigma_{V,V},\quad
\frosp(S_{1,2},b)=\sigma_{1,V}
\]
(note that $\sigma_{1,1}=0$). But for any $u,v,w\in V$:
\begin{equation}\label{eq:sigmauvB12}
\sigma_{u,v}:\begin{cases} 1\mapsto 0,&\\
 w\mapsto -b(u,w)v-b(v,w)u=-\gamma_{u,v}(w),& \end{cases}
\end{equation}
while
\begin{equation}\label{eq:sigma1uB12}
\sigma_{1,u}:\begin{cases} 1\mapsto 2u=-u\ (\charac k=3),&\\
  v\mapsto -b(u,v)1=-\langle u\vert v\rangle 1.& \end{cases}
\end{equation}

Also, $[\sigma_{1,u},\sigma_{1,v}]=-\sigma_{u,v}$ and
$[\sigma_{u,v},\sigma_{1,w}]=\sigma_{1,\sigma_{u,v}(w)}$ for any
$u,v,w\in V$. Thus, $\sigma_{V,V}$ can be identified with the three
dimensional symplectic Lie algebra $\frsp(V)$, $\sigma_{1,V}$ with
$V$ ($\sigma_{1,u}\leftrightarrow u$), and therefore
$\frosp(S_{1,2},b)$ is isomorphic to the Lie superalgebra
\begin{equation}\label{eq:b01}
\frb_{0,1}=\frsp(V)\oplus V,
\end{equation}
with even part $\frsp(V)$, odd part $V$, and multiplication
determined by:
\begin{itemize}
\item the Lie algebra $\frsp(V)$ is a subalgebra,
\smallskip

\item $[\gamma, v]=\gamma(v)$ for any $\gamma\in \frsp(V)$ and $v\in
V$,
\smallskip

\item $[u,v]=\gamma_{u,v}$ for any $u,v\in V$.
\end{itemize}

\bigskip

Now, consider the superalgebra $S_{4,2}=\overline{B(4,2)}$ over a
field $k$ of characteristic $3$, as described in Example
\ref{ex:B42bar}. Then \cite[Lemma 5.7]{EldOkuSuperCompo} shows that
the Lie superalgebra $\frosp(S_{4,2},b)$ is isomorphic to the Lie
superalgebra
\begin{equation}\label{eq:d21}
\frd_{2,1}=\bigl(\frsp(V)\oplus\frsp(V)\oplus\frsp(V)\bigr)\oplus
(V\otimes V\otimes V),
\end{equation}
with even part $\frsp(V)\oplus\frsp(V)\oplus\frsp(V)$, odd part
$V\otimes V\otimes V$, and multiplication determined by:
\begin{itemize}
\item the Lie algebra $\frsp(V)\oplus\frsp(V)\oplus\frsp(V)$ is a
subalgebra,
\smallskip

\item for any
$\gamma_1,\gamma_2,\gamma_3\in \frsp(V)$ and $u_1,u_2,u_3\in V$,
\begin{multline*}
[(\gamma_1,\gamma_2,\gamma_3),u_1\otimes u_2\otimes u_3]\\
=\gamma_1(u_1)\otimes u_2\otimes u_3+u_1\otimes \gamma_2(u_2)\otimes
u_3+u_1\otimes u_2\otimes \gamma(u_3),
\end{multline*}

\item for any $u_1,u_2,u_3,v_1,v_2,v_3\in V$,
\[
\begin{split}
[u_1&\otimes u_2\otimes u_3,v_1\otimes v_2\otimes v_3]\\
 &=-\bigl(\langle u_2\vert v_2\rangle\langle u_3\vert v_3\rangle
     \gamma_{u_1,v_1},
     \langle u_1\vert v_1\rangle\langle u_3\vert v_3\rangle
     \gamma_{u_2,v_2},
     \langle u_1\vert v_1\rangle\langle u_2\vert v_2\rangle
     \gamma_{u_3,v_3}\bigr)\\
 &=-\sum_{i=1}^3\bigl(\prod_{j\ne i}\langle u_j\vert
 v_j\rangle\bigr)\nu_i(\gamma_{u_i,v_i}),
\end{split}
\]
where $\nu_i: \frsp(V)\rightarrow \frsp(V)^3$ denotes the inclusion
on the $i^{\text{\small th}}$-component.
\end{itemize}

(A word of caution is needed here: the operators $\sigma_{a,b}$ in
\cite{EldOkuSuperCompo} are changed in sign with respect to ours.)

\smallskip

Moreover, the action of $\frd_{2,1}$ on $S_{4,2}=(V\otimes V)\oplus
V$ is given by the isomorphism:
\begin{equation}\label{eq:rhod21}
\rho:\frd_{2,1}\rightarrow \frosp(S_{4,2},b)
\end{equation}
such that
\[
\begin{split}
\rho\bigl((\gamma_1,\gamma_2,\gamma_3)\bigr)(v_1\otimes v_2 + v_3)&=
  \bigl(\gamma_1(v_1)\otimes v_2+v_1\otimes
  \gamma_2(v_2)\bigr)+\gamma_3(v_3),\\[8pt]
\rho\bigl(u_1\otimes u_2\otimes u_3)(v_1\otimes v_2 + v_3)&=
  -\sigma_{u_1\otimes u_2,u_3}(v_1\otimes v_2 + v_3)\\
  &=\langle u_3\vert v_3\rangle u_1\otimes u_2 -
    \langle u_1\vert v_1\rangle\langle u_2\vert v_2\rangle u_3,
\end{split}
\]
for any $\gamma_1,\gamma_2,\gamma_3\in \frsp(V)$ and
$u_1,u_2,u_3,v_1,v_2,v_3\in V$.

Note that
$\sigma_{u_3,v_3}(w_3)=-b(u_3,w_3)v_3-b(v_3,w_3)u_3=-\gamma_{u_3,v_3}(w_3)$,
so (see \cite[(2.8)]{ElduqueMagicII})
\begin{equation}\label{eq:sigmasS42}
\begin{split}
\sigma_{u_1\otimes u_2,v_1\otimes v_2}
 &=-\rho\bigl(\langle u_2\vert v_2\rangle \gamma_{u_1,v_1},
  \langle u_1\vert v_1\rangle\gamma_{u_2,v_2},0)\bigr),\\
\sigma_{u_3,v_3}&=-\rho\bigl((0,0,\gamma_{u_3,v_3})\bigr),\\
\sigma_{u_1\otimes u_2,u_3}&=-\rho\bigl(u_1\otimes u_2\otimes
u_3\bigr),
\end{split}
\end{equation}
for any $u_1,u_2,u_3,v_1,v_2,v_3\in V$.

\smallskip

Consider the natural order $3$ automorphism $\theta$ of $\frd_{2,1}$
such that
\begin{equation}\label{eq:theta}
\left\{\begin{aligned}
 &\theta(\gamma_1,\gamma_2,\gamma_3)=(\gamma_3,\gamma_1,\gamma_2),\\
 &\theta(u_1\otimes u_2\otimes u_3)=u_3\otimes u_1\otimes u_2,
\end{aligned}\right.
\end{equation}
for any $\gamma_1,\gamma_2,\gamma_3\in \frsp(V)$ and $u_1,u_2,u_3\in
V$.

Then (compare to \cite[Proposition 2.12]{ElduqueMagicII}):

\begin{proposition}\label{pr:trialityd21}
For any homogeneous elements $f\in \frd_{2,1}$ and $x,y\in S_{4,2}$:
\[
\rho(f)(x\bullet y)=\rho(\theta^{-1}(f))(x)\bullet y+
 (-1)^{\lvert f\rvert\lvert x\rvert}x\bullet
 \rho(\theta^{-2}(f))(y).
\]
\end{proposition}
\begin{proof}
It is enough to prove this for generators of $\frd_{2,1}$ and a
spanning set of $S_{4,2}$, and hence for $f=u_1\otimes u_2\otimes
u_3$,  $x=v_1\otimes v_2$ or $x=v_3$, and $y=w_1\otimes w_2$ or
$y=w_3$, where $u_i,v_i,w_i\in V$, $i=1,2,3$. This is
straightforward. For instance,
\[
\begin{split}
\rho(u_1\otimes u_2\otimes u_3)
  \bigl(v_3\bullet(w_1\otimes w_2)\bigr)&=
   \rho(u_1\otimes u_2\otimes u_3)\bigl(\langle w_2\vert v_3\rangle
   w_1\bigr)\\
  &=\langle w_2\vert v_3\rangle\langle u_3\vert w_1\rangle
  u_1\otimes u_2,\\[3pt]
\rho\bigl(\theta^{-1}(u_1\otimes u_2\otimes u_3)\bigr)(v_3)\bullet
      (w_1\otimes w_2)&=
   \rho(u_2\otimes u_3\otimes u_1)(v_3)\bullet (w_1\otimes w_2)\\
   &=\langle u_1\vert v_3\rangle (u_2\otimes u_3)\bullet (w_1\otimes
   w_2)\\
   &=\langle u_1\vert v_3\rangle \langle u_3\vert w_1\rangle
   w_2\otimes u_2\\[3pt]
-v_3\bullet\rho\bigl(\theta^{-2}(u_1\otimes u_2\otimes
           u_3)\bigr)(w_1\otimes w_2)&=
  -v_3\bullet\bigl(\rho(u_3\otimes u_1\otimes u_2)(w_1\otimes w_2)\\
  &=v_3\bullet\bigl(\langle u_3\vert w_1\rangle\langle u_1\vert
  w_2\rangle u_2\\
  &=-\langle u_3\vert w_1\rangle\langle u_1\vert w_2\rangle
  v_3\otimes u_2,
\end{split}
\]
and now, since $\langle u\vert v\rangle w+\langle v\vert w\rangle u+
\langle w\vert u\rangle v=0$ for any $u,v,w\in V$ (as any trilinear
alternating form on a two dimensional vector space is trivially
zero), one gets
\[
\begin{split}
 \rho(u_1\otimes u_2\otimes u_3)\bigl(v_3\bullet(w_1\otimes
 w_2)\bigr)
  &=\rho\bigl(\theta^{-1}(u_1\otimes u_2\otimes u_3)\bigr)(v_3)\bullet
      (w_1\otimes w_2) \\
      &\quad - v_3\bullet\rho\bigl(\theta^{-2}(u_1\otimes u_2\otimes
           u_3)\bigr)(w_1\otimes w_2),
\end{split}
\]
as required.
\end{proof}

\medskip

Write $\rho_i=\rho\circ\theta^{-i}$, then Theorem
\ref{th:triality1224} and Proposition \ref{pr:trialityd21} give, as
in \cite[Corollary 2.13]{ElduqueMagicII}:

\begin{corollary}
$\tri(S_{4,2})=\bigl\{\bigl(\rho_0(f),\rho_1(f),\rho_2(f)\bigr):
f\in\frd_{2,1}\bigr\}$.
\end{corollary}

This Corollary allows us to identify $\tri(S_{4,2})$ to $\frd_{2,1}$
by means of $f\leftrightarrow
\bigl(\rho_0(f),\rho_1(f),\rho_2(f)\bigr)$.

\begin{remark}
The formulas for the bracket in $\frd_{2,1}$ and for the
representations $\rho_i$ of $\frd_{2,1}$ ($i=0,1,2$) are exactly the
ones that appear in \cite[\S 2]{ElduqueMagicII} for $\frd_4$, but
with the last copy of $V$ suppressed.
\end{remark}

\medskip

Denote by $V_i$ ($i=1,2,3$) the module $V$ for
$\frsp(V)^3=\frsp(V)\oplus\frsp(V)\oplus\frsp(V)$ on which only the
$i^{\text{\small th}}$ component of $\frsp(V)^3$ acts:
$(s_1,s_2,s_3).v_i=s_i(v_i)$, for any $s_j\in\frsp(V)$, $j=1,2,3$,
and $v_i\in V_i$. Also denote by $\iota_i(S_{4,2})$ the
$\frd_{2,1}$-module associated to the representation $\rho_i$. Then,
as modules for $\frsp(V)^3$ (compare to
\cite[(2.14)]{ElduqueMagicII}):
\begin{equation}\label{eq:iotaS42}
\left\{
\begin{aligned}
\iota_0(S_{4,2})&=(V_1\otimes V_2)\oplus V_3,\\
\iota_1(S_{4,2})&=(V_2\otimes V_3)\oplus V_1,\\
\iota_2(S_{4,2})&=(V_3\otimes V_1)\oplus V_2\,\bigl(\simeq
(V_1\otimes V_3)\oplus V_2\bigr).
\end{aligned}\right.
\end{equation}

The multiplication $\bullet$ on $S_{4,2}$ (see Example
\ref{ex:B42bar}) becomes the bilinear $\tri(S_{4,2})$-invariant map
\[
\iota_0(S_{4,2})\times \iota_1(S_{4,2})\rightarrow \iota_2(S_{4,2})
\]
given by,
\begin{equation}\label{eq:bulletiotaS42}
\begin{split}
\bigl((V_1\otimes V_2)\oplus V_3\bigr)\times \bigl(V_2\otimes
V_3)\oplus V_1\bigr)&\rightarrow (V_3\otimes V_1)\oplus V_2\\
\bigl(u_1\otimes u_2 +u_3),(v_2\otimes v_3+v_1\bigr)&\mapsto
 \bigl(\langle u_2\vert v_2\rangle v_3\otimes u_1 -u_3\otimes v_1\bigr)\\
 &\null\qquad\qquad\qquad
  -\bigl(\langle u_1\vert v_1\rangle u_2+\langle u_3\vert v_3\rangle
  v_2\bigr)
\end{split}
\end{equation}
for any $u_i,v_i\in V_i$, $i=1,2,3$; and cyclically. Note that it
consists of contractions for repeated indices.

This is exactly the multiplication in \cite[3.1]{ElduqueMagicII} if
$V_4$ is ignored there.

\bigskip

\section{The extended Freudenthal's Magic Square}

Let $(S,\bullet,q)$ and $(S',*,q')$ be two symmetric composition
superalgebras and define $\frg=\frg(S,S')$ to be the
$\bZ_2\times\bZ_2$-graded anticommutative superalgebra such that
\[
\begin{split}
&\frgoo=\tri(S,\bullet,q)\oplus\tri(S',*,q'),\\
&\frgunoo=\frgouno=\frgunouno=S\otimes S'.
\end{split}
\]

For any $x\in S$ and $x'\in S'$, denote by $\iota_i(x\otimes x')$
the element $x\otimes x'$ in $\frgunoo$ (respectively $\frgouno$,
$\frgunouno$) if $i=0$ (respectively, $i=1,2$). Thus
\begin{equation}\label{eq:gSSprime}
\frg=\frg(S,S')=\bigl(\tri(S,\bullet,q)\oplus\tri(S',*,q')\bigr)
\oplus \bigl(\oplus_{i=0}^2\iota_i(S\otimes S')\bigr)\,.
\end{equation}
Then $\frg$ is a superalgebra where, for $i=\bar 0,\bar 1$,
\[
\frg_i=(\frgoo)_i\oplus
(\frgunoo)_i\oplus(\frgouno)_i\oplus(\frgunouno)_i,
\]
with
\[
\begin{split}
&(\frgoo)_i=\tri(S)_i\oplus\tri(S')_i,\\
&\iota_j(S\otimes S')_i=\iota_j(S\subo\otimes S'_i)\oplus
\iota_j(S\subuno\otimes S'_{\bar 1-i}),
\end{split}
\]
for $j=0,1,2$.

 The
superanticommutative multiplication on $\frg$ is defined by means
of:
\begin{itemize}
\item $\frgoo$ is a Lie subsuperalgebra of $\frg$,
\smallskip

\item $[(d_0,d_1,d_2),\iota_i(x\otimes
 x')]=\iota_i\bigl(d_i(x)\otimes x'\bigr)$,
\smallskip

\item
 $[(d_0',d_1',d_2'),\iota_i(x\otimes
 x')]=(-1)^{\lvert d_i'\rvert\lvert x\rvert}\iota_i\bigl(x\otimes d_i'(x')\bigr)$,
\smallskip

\item $[\iota_i(x\otimes x'),\iota_{i+1}(y\otimes y')]=(-1)^{\lvert
x'\rvert\lvert y\rvert}
 \iota_{i+2}\bigl((x\bullet y)\otimes (x'*y')\bigr)$ (indices modulo
 $3$),
\smallskip

\item $[\iota_i(x\otimes x'),\iota_i(y\otimes y')]=
 (-1)^{\lvert x\rvert\lvert x'\rvert+\lvert x\rvert\lvert y'\rvert +
 \lvert y\rvert\lvert y'\rvert}
 b'(x',y')\theta^i(t_{x,y})$ \newline \null\qquad\qquad\qquad\qquad $+
 (-1)^{\lvert y\rvert\lvert x'\rvert}
 b(x,y)\theta'^i(t'_{x',y'})$,

\end{itemize}
for any $i=0,1,2$ and homogeneous $x,y\in S$, $x',y'\in S'$,
$(d_0,d_1,d_2)\in\tri(S)$, and $(d_0',d_1',d_2')\in\tri(S')$. Here
$\theta$ denotes the natural automorphism
$\theta:(d_0,d_1,d_2)\mapsto (d_2,d_0,d_1)$ in $\tri(S)$, $\theta'$
the analogous automorphism of $\tri(S')$, and
\begin{equation}\label{eq:txy}
t_{x,y}=\bigl(\sigma_{x,y},\tfrac{1}{2}b(x,y)1-r_xl_y,\tfrac{1}{2}b(x,y)1-l_xr_y\bigr)
\end{equation}
(with $l_x(y)=x\bullet y$, $r_x(y)=(-1)^{\lvert x\rvert\lvert y
\rvert}y\bullet x$), while $t'_{x',y'}$ is the analogous triple in
$\tri(S')$.

\medskip

Just superizing the arguments in \cite[Theorem 3.1]{ElduqueMagicI};
that is, taking into account the parity signs, one gets:

\begin{theorem}
With this multiplication, $\frg=\frg(S,S')$ is a Lie superalgebra.
\end{theorem}

\medskip

With the same proof as in the ungraded case \cite[Theorem
12.2]{ElduqueNewLook} the following result is obtained:

\begin{proposition}\label{pr:SSphi}
Let $S$ and $S'$ be two symmetric composition superalgebras, and let
$\varphi$ be an automorphism of $S$ of order $3$, then the Lie
superalgebras $\frg(S,S')$ and $\frg(S_\varphi,S')$ are isomorphic.
\end{proposition}

(The superalgebra $S_\varphi$ was defined in \eqref{eq:Sphi}.)

\medskip

Therefore, over fields of characteristic $3$, there is no need to
deal with the symmetric composition superalgebras
$\overline{B(1,2)}_\lambda$, but just with
$S_{1,2}=\overline{B(1,2)}$ (which is three dimensional) and with
$S_{4,2}=\overline{B(4,2)}$ (whose dimension is $6$). Freudenthal's
Magic Square thus extends over these fields to the larger square
$$
\vbox{\offinterlineskip
 \halign{$#$\qquad&\hfil$#$\quad\hfil&%
 \vreglon #%
 &\hfil\quad$#$\quad\hfil&\hfil\quad$#$\quad\hfil
 &\hfil\quad$#$\quad\hfil&\hfil\quad$#$\quad\hfil&%
 \vrule  depth 4pt width .5pt #%
 &\hfil\quad$#$\quad\hfil&\hfil\quad$#$\ \hfil\cr
 &&\omit&\multispan7{\hfil$\dim S$\qquad\hfil}\cr
 \noalign{\medskip}
 \bigstrut &&width 0pt&1&2&4&8&\omit\vrule height 8pt depth 4pt width .5pt&3&6\cr
 &&\multispan8{\hreglonfill}\cr
 &1&&A_1&\tilde A_2&C_3&F_4&&\frg(S_1,S_{1,2})&\frg(S_1,S_{4,2})\cr
 \bigstrut&2&& &\omit$\tilde A_2\oplus \tilde A_2$&\tilde A_5&\tilde E_6&&\frg(S_2,S_{1,2})&\frg(S_2,S_{4,2})\cr
 \bigstrut\smash{\lower 6pt\hbox{$\dim S'$}}&4&& & &D_6&E_7&&\frg(S_4,S_{1,2})&\frg(S_4,S_{4,2})\cr
 \bigstrut&8&& & & &E_8&&\frg(S_8,S_{1,2})&\frg(S_8,S_{4,2})\cr
 &\multispan9{\hregletafill}\cr
 \bigstrut& 3&& & & & & & \frg(S_{1,2},S_{1,2})&\frg(S_{1,2},S_{4,2})\cr
 \bigstrut& 6&& & & & & & & \frg(S_{4,2},S_{4,2})\cr}}
$$
\bigskip

\noindent (see Table \ref{ta:superFreudenthal}, where the dimensions
of the even and odd part of the superalgebras involved, which are
easily computed from the definitions, are displayed), where $S_n$
denotes a symmetric composition algebra of dimension $n$ ($n=1,2,4$
or $8$).

\bigskip

The purpose for the remaining part of the paper is the description
of the Lie superalgebras $\frg(S,S_{1,2})$ and $\frg(S,S_{4,2})$ in
the split case, that is, in case the (super)algebra $S$ contains
nontrivial idempotents (for instance, over an algebraically closed
field). Proposition \ref{pr:SSphi}, or \cite[Corollary
12.1]{ElduqueNewLook} shows that it is enough to take then for $S_n$
the unique split para-Hurwitz algebras (dimensions $1$, $2$, $4$ or
$8$), together with $S_{1,2}$ and $S_{4,2}$.

\smallskip

This description will be given in terms of contragredient Lie
superalgebras.

\bigskip

\section{Contragredient Lie superalgebras}

Let $A=\bigl(a_{ij}\bigr)$ be a square matrix of size $n$ over a
ground field $k$, and let $(\frh,\Pi,\Pi^{\vee})$ be a
\emph{realization} of $A$, as defined in \cite[\S 1.1]{KacBook}.
That is,
\begin{itemize}
\item $\Pi=\{\alpha_1,\ldots,\alpha_n\}$ is a linearly independent
set in $\frh^*$ (the dual of the vector space $\frh$),

\item $\Pi^{\vee}=\{h_1,\ldots,h_n\}$ is a linearly independent set
in $\frh$,

\item $\alpha_j(h_i)=a_{ij}$ for any $i,j=1,\ldots,n$,

\item $\dim\frh=2n-\rank A$.
\end{itemize}

As in \cite{KacLieSuper}, for any subset $\tau\subseteq \{
1,\ldots,n\}$, consider the local Lie superalgebra
\[
\hat\frg(A,\tau)=\frg_{-1}\oplus\frg_0\oplus \frg_1,
\]
with $\frg_0=\frh$, $\frg_1=ke_1\oplus\cdots\oplus ke_n$,
$\frg_{-1}=kf_1\oplus\cdots \oplus kf_n$, and bracket given by:
\begin{equation}\label{eq:realization}
[e_i,f_j]=\delta_{ij}h_i,\quad
 [h,h']=0,\quad
 [h,e_i]=\alpha_i(h)e_i,\quad
 [h,f_i]=-\alpha_i(h)f_i,
\end{equation}
for any $i,j=1,\ldots,n$ and $h,h'\in\frh$, where $\frg_0$ is even
and $e_i,f_i$ are even if and only if $i\not\in\tau$.

Note that changing $f_i$ by $cf_i$, $0\ne c\in k$, the
$i^{\text{\small th}}$ row of $A$ is multiplied by $c$. Hence, if
$a_{ii}\ne 0$, then $a_{ii}$ can be taken to be $2$, as it is
customary.

Then (see \cite[1.2.2]{KacLieSuper}) there exists a minimal
$\bZ$-graded Lie superalgebra $\frg(A,\tau)$ with local part
$\hat\frg(A,\tau)$. To define $\frg(A,\tau)$ one first considers
$\tilde\frg(A,\tau)$, the free Lie superalgebra generated by $\frh$
and $\{e_1,\ldots,e_n,f_1,\ldots,f_n\}$, subject to the relations in
\eqref{eq:realization}, endowed with the $\bZ$-grading induced by
the grading in $\hat\frg(A,\tau)$ (note that the relations are
homogeneous), and then considers the largest homogeneous ideal
$\fri(A)$ that intersects $\frh$ trivially. Then
\[
\frg(A,\tau)=\tilde\frg(A,\tau)/\fri(A).
\]

\begin{remark}\null\quad\null
\begin{itemize}
\item As in \cite[Proposition 1.6]{KacBook}, with
$\frg'(A,\tau)=[\frg(A,\tau),\frg(A,\tau)]$, the center is given by
\[
\frc=Z\bigl(\frg(A,\tau)\bigr)=Z\bigl(\frg'(A,\tau)\bigr)=
 \{ h\in \frh : \alpha_i(h)=0\ \forall i=1,\ldots,n\}.
\]
Note that $\dim\frc=n-\rank A$.

\item Let $Q$ be the free abelian group on generators
$\epsilon_1,\ldots,\epsilon_n$. Then $\frg(A,\tau)$ is $Q$-graded by
assigning $\deg e_i=\epsilon_i=-\deg f_i$, $i=1,\ldots,n$, and $\deg
h=0$ for any $h\in\frh$. Now $\tilde\frg(A,\tau)$ is $Q$-graded too,
because the relation in \eqref{eq:realization} are $Q$-homogeneous.
For any $m\in\bZ$, the $m^{\text{\small th}}$ homogeneous component
of $\tilde\frg(A,\tau)$ equals
\[
\bigoplus_{\substack{q=m_1\epsilon_1+\cdots+m_n\epsilon_n\\
 m_1+\cdots+m_n=m}} \tilde\frg(A,\tau)_q,
\]
and the ideal $\fri(A)$ is $Q$-homogeneous too, since $\fri(A)$ is
contained in $\oplus_{0\ne q\in Q}\pi_q\bigl(\fri(A)\bigr)$ ($\pi_q$
denotes the projection onto the $q^{\text{\small th}}$-component),
which is a $Q$-homogeneous ideal that intersects trivially $\frh$.
Hence, by maximality, $\fri(A)=\oplus_{0\ne q\in
Q}\pi_q\bigl(\fri(A)\bigr)$.

\end{itemize}
\end{remark}

\medskip

The following notation will be used:
\begin{itemize}
\item[$\diamond$] $\frg(A,\tau)$ will be called the \emph{contragredient Lie
superalgebra} with Cartan matrix $A$.

\item[$\diamond$] $\frg'(A,\tau)$ will be called the \emph{derived
contragredient Lie superalgebra} with Cartan matrix $A$.

\item[$\diamond$] $\frg(A,\tau)/\frc$ will be called the \emph{centerless
contragredient Lie superalgebra} with Cartan matrix $A$.

\item[$\diamond$] $\frg'(A,\tau)/\frc$ will be called the \emph{centerless
derived contragredient Lie superalgebra} with Cartan matrix $A$.

\end{itemize}

A bit of caution is needed here, in \cite{VeisfeilerKac} it is the
Lie algebra $\frg'(A,\tau)$ the one that is called the
contragredient Lie algebra.

\medskip

Note that, as in \cite[Proposition 1.7]{KacBook}, $\frg(A,\tau)$ has
no $Q$-homogeneous ideal if and only if $\det A\ne 0$ and for any
$i,j$ there exist indices $i_1,\ldots, i_s$ such that
$a_{ii_1}a_{i_1i_2}\cdots a_{i_sj}\ne 0$. If only this last
condition holds, then any $Q$-homogeneous ideal either contains
$\frg'(A,\tau)$ or is contained in the center. Also note that if
$\det A\ne 0$, then $\frc=0$, $\frg'(A,\tau)=\frg(A,\tau)$ and the
four Lie superalgebras considered above coincide.

\begin{lemma}\label{le:jA} Under the conditions above,
let $\frj(A)$ be the largest $\bZ$-homoge\-ne\-ous ideal of
$\tilde\frg(A,\tau)$ which intersects trivially
$\frg_{-1}\oplus\frg_1$. Then $\frj(A)=\fri(A)\oplus\frc$.
\end{lemma}
\begin{proof}
It is clear that $\fri(A)\cap\frg_{\pm 1}=0$, otherwise $\fri(A)$
would contain some $0\ne h\in\frh$ (because of $[e_i,f_i]=h_i$ for
any $i$). Hence $\fri(A)$ is the largest homogeneous ideal
intersecting trivially $\frg_{-1}\oplus\frg_0\oplus\frg_1$. Thus
$\bigl(\fri(A)\oplus\frc\bigr)\cap\bigl(\frg_{-1}\oplus\frg_1\bigr)=0$,
so $\fri(A)\oplus\frc\subseteq \frj(A)$.

Conversely, $\tilde\frg(A,\tau)=\oplus_{i\in\bZ}\tilde\frg_i$, with
$\tilde\frg_i=\frg_i$ for $i=-1,0,1$, and for any
$h\in\frj(A)\cap\tilde\frg_0=\frj(A)\cap \frh$,
$[h,e_i]=\alpha_i(h)e_i\in\frj(A)\cap\frg_1=0$ for any $i$. Thus
$\alpha_i(h)=0$ for any $i$ and $h\in\frc$. Therefore
$\frj(A)\subseteq \frc\oplus\bigl(\oplus_{i\ne
-1,0,1}\frj(A)\cap\tilde\frg_i\bigr)$. But $\oplus_{i\ne
-1,0,1}\frj(A)\cap\tilde\frg_i$ is an ideal of $\tilde\frg(A,\tau)$,
since it is closed under the action of $e_i,f_i$, $i=1,\ldots,n$,
and $\frh$ (note that $[f_i,\frj(A)\cap\tilde\frg_2]\subseteq
\frj(A)\cap \frg_1=0$). Hence $\oplus_{i\ne
-1,0,1}\frj(A)\cap\tilde\frg_i\subseteq \fri(A)$, and
$\frj(A)\subseteq\fri(A)\oplus\frc$.
\end{proof}

\smallskip

Let us consider the natural characterizations of the contragredient
Lie superalgebras, that will be used in the next section. The
previous notations and assumptions will be kept in the following
results.

\begin{theorem}\label{th:gAtau}
Let $\frg$ be a $\bZ$-graded Lie superalgebra generated by $\frg_0$,
which is contained in $\frg\subo$, and by elements $e_i\in\frg_1$
and $f_i\in\frg_{-1}$, $i=1,\ldots,n$, which are even (respectively
odd) if and only if $i\not\in\tau$ (resp. $i\in\tau$). Assume that
there are $\alpha_1,\ldots,\alpha_n\in \frg_0^*$ and
$h_1,\ldots,h_n\in \frg_0$ such that
\[
\bigl(\frg_0,\{\alpha_1,\ldots,\alpha_n\},\{h_1,\ldots,h_n\}\bigr)
\]
is a realization of $A$ and that relations \eqref{eq:realization}
hold. If any nonzero homogeneous ideal of $\frg$ intersects $\frg_0$
nontrivially, then $\frg$ is isomorphic to the contragredient Lie
superalgebra $\frg(A,\tau)$.
\end{theorem}
\begin{proof}
This is clear from the definitions, which give an epimorphism
$\phi:\frg(A,\tau)\rightarrow \frg$, such that
$\ker\phi\cap\frg(A,\tau)_0=0$, so $\ker\phi=0$.
\end{proof}

\medskip

\begin{theorem}\label{th:centerlessgAtau}
Let $\frg$ be a $\bZ$-graded Lie superalgebra such that:
\begin{romanenumerate}
\item $\frg_0=\bar\frh$ is contained in $\frg\subo$, it is abelian
and its dimension is $n$.

\item There are linearly independent elements
$\bar\alpha_1,\ldots,\bar\alpha_n\in\bar\frh^*$ and elements $\bar
h_1,\ldots,\bar h_n\in\bar\frh$ such that $\bar\alpha_j(\bar
h_i)=a_{ij}$ for any $i,j$.

\item There are elements $\bar e_1,\ldots,\bar e_n\in\frg_1$ and
$\bar f_1,\ldots,\bar f_n\in\frg_{-1}$, where $\bar e_i,\bar f_i$
are even (resp. odd) if $i\not\in\tau$ (resp. $i\in\tau$), and such
that the relations \eqref{eq:realization} are satisfied (with $e_i$
replaced by $\bar e_i$, ...), and $\frg$ is generated by $\bar
e_1,\ldots,\bar e_n,\bar f_1,\ldots,\bar f_n$ and $\bar\frh$.

\item Any nonzero homogeneous ideal of $\frg$ intersects
$\frg_{-1}\oplus\frg_1$ nontrivially.

\end{romanenumerate}

Then $\frg$ is isomorphic to the centerless contragredient Lie
superalgebra $\frg(A,\tau)/\frc$.
\end{theorem}
\begin{proof}
The elements $\alpha_1,\ldots,\alpha_n\in\frh^*$
($\frh=\frg(A,\tau)_0$) define linearly independent elements
$\hat\alpha_i\in (\frh/\frc)^*$, because $\frc=\{h\in\frh:
\alpha_i(h)=0\ \forall i\}$. Hence there is a linear bijection
$\varphi^*:\bar\frh^*\rightarrow (\frh/\frc)^*$ such that
$\varphi^*(\bar\alpha_i)=\hat\alpha_i$ for any $i$. Besides,
$\varphi^*$ induces a linear bijection $\varphi:
\frh/\frc\rightarrow \bar\frh$, such that
$\bar\alpha_j\bigl(\varphi(h+\frc)\bigr)=\hat\alpha_j(h+\frc)=\alpha_j(h)$
for any $j$ and $h\in\frh$. By composing with the natural
projection, a surjective linear map is obtained
$\phi:\frh\rightarrow \bar\frh$ such that
$\bar\alpha_j\bigl(\phi(h)\bigr)=\alpha_j(h)$ for any $j$ and $h\in
\frh$.

Now, by definition of $\tilde\frg(A,\tau)$, $\phi$ extends to a
surjective homomorphism of Lie superalgebras
$\tilde\phi:\tilde\frg(A,\tau)\rightarrow \frg$ which takes $e_i$
(respectively $f_i$) to $\bar e_i$ (respectively $\bar f_i$), for
any $i$, and $\tilde\phi$ is homogeneous. Condition (iv) shows that
$\tilde\phi\bigl(\frj(A)\bigr)=0$, so $\tilde \phi$ induces a
surjective homomorphism
$\hat\phi:\tilde\frg(A,\tau)/\frj(A)\rightarrow \frg$. But
$\hat\phi$ is bijective in degrees $0$, $1$ and $-1$ (the elements
$\bar e_1,\ldots,\bar e_n$ are linearly independent, as they are
eigenvectors of eigenvalues $\bar\alpha_1,\ldots,\bar\alpha_n$ for
the action of $\bar\frh$). By definition of $\frj(A)$, one concludes
that $\ker\hat\phi=0$, so $\frg$ is isomorphic to
$\tilde\frg(A,\tau)/\frj(A)$, which is isomorphic to
$\frg(A,\tau)/\frc$ by Lemma \ref{le:jA}.
\end{proof}

\medskip

\begin{theorem}\label{th:derivedgAtau}
Let $\frg$ be a $Q$-graded Lie superalgebra generated by nonzero
elements $\bar e_i\in\frg_{\epsilon_i}$, $\bar
f_i\in\frg_{-\epsilon_i}$, $i=1,\ldots,n$, where $\bar e_i,\bar f_i$
are even (resp. odd) if $i\not\in\tau$ (resp. $i\in\tau$). Assume
that:
\begin{romanenumerate}
\item The elements $\bar h_i=[\bar e_i,\bar f_i]$, $i=1,\ldots,n$,
are linearly independent.

\item $[\bar h_i,\bar e_j]=a_{ij}\bar e_j$, $[\bar h_i,\bar
f_j]=-a_{ij}\bar f_j$, $[\bar e_i,\bar f_j]=\delta_{ij}\bar h_i$,
and $[\bar h_i,\bar h_j]=0$, for any $i,j=1,\ldots,n$.

\item Any nonzero $Q$-homogeneous ideal of $\frg$ intersects
nontrivially $\frg_{-1}\oplus\frg_{0}\oplus\frg_1$, where $\frg_m=
\oplus \{\tilde\frg_q :
 q=m_1\epsilon_1+\cdots+m_n\epsilon_n,\,
 m_1+\cdots+m_n=m\}$ for any $m\in\bZ$.

\end{romanenumerate}

Then $\frg$ is isomorphic to the derived contragredient Lie
superalgebra $\frg'(A,\tau)$.
\end{theorem}
\begin{proof}
Since $\frg$ is generated by the elements $\bar e_1,\ldots,\bar
e_n,\bar f_1,\ldots,\bar f_n$, it follows that $\frg_0=k\bar
h_1\oplus\cdots\oplus k\bar h_n$. Take an even vector space
$\hat\frg_0$ of dimension $n-\rank A$ and elements $\bar
\alpha_1,\ldots,\bar\alpha_n\in (\frg_0\oplus\hat\frg_0)^*$ such
that
\[
\left\{\begin{aligned} &\bar\alpha_j(\bar h_i)=a_{ij}\ \text{for any
$i,j$,}\\
 &\text{$\bar\alpha_1,\ldots,\bar\alpha_n$ are linearly independent.}
 \end{aligned}\right.
\]
Then
$\bigl(\frg_0\oplus\hat\frg_0,\{\bar\alpha_1,\ldots,\bar\alpha_n\},\{\bar
h_1,\ldots,\bar h_n\}\bigr)$ is a realization of $A$, and
$\bar\frg=\frg\oplus\hat\frg_0$ is a $Q$-graded Lie superalgebra ,
with $\bar\frg_0=\frg_0\oplus\hat\frg_0$,
$\bar\frg_\epsilon=\frg_\epsilon$ for any $0\ne\epsilon\in Q$, and
where $[\bar\frg_0,\bar\frg_0]=0$, $\frg$ is an ideal of $\bar\frg$,
and for any $h\in \hat\frg_0$ and $x\in\frg_\epsilon$
($0\ne\epsilon=m_1\epsilon_1+\cdots+m_n\epsilon_n\in Q$),
$[h,x]=(m_1\bar\alpha_1+\cdots+m_n\bar\alpha_n)(h)x$.

By Theorem \ref{th:gAtau}, $\bar\frg$ is isomorphic to the
contragredient Lie superalgebra $\frg(A,\tau)$, and hence
$\frg=[\bar\frg,\bar\frg]$ is isomorphic to $\frg'(A,\tau)$.
\end{proof}

\medskip

\begin{theorem}\label{th:centerlessderivedgAtau}
Let $\frg$ be a $Q$-graded Lie superalgebra generated by nonzero
elements $\bar e_i\in\frg_{\epsilon_i}$, $\bar f_i\in
\frg_{-\epsilon_i}$, $i=1,\ldots,n$, where $\bar e_i,\bar f_i$ are
even (resp. odd) if $i\not\in\tau$ (resp. $i\in\tau$). Assume that:
\begin{romanenumerate}
\item If $\bar h_i=[\bar e_i,\bar f_i]$, $i=1,\ldots,n$, then $[\bar
h_i,\bar e_j]=a_{ij}\bar e_j$, $[\bar h_i,\bar f_j]=-a_{ij}\bar
f_j$, $[\bar h_i,\bar h_j]=0$, and $[\bar e_i,\bar
f_j]=\delta_{ij}\bar h_i$, for any $i,j$.

\item Any nonzero $Q$-homogeneous ideal intersects nontrivially
$\frg_{-1}\oplus \frg_1$, where $\frg_m= \oplus \{\tilde\frg_q :
q=m_1\epsilon_1+\cdots+m_n\epsilon_n,\,
 m_1+\cdots+m_n=m\}$ for any $m\in\bZ$.

\end{romanenumerate}

Then $\frg$ is isomorphic to the centerless derived contragredient
Lie superalgebra $\frg'(A,\tau)/\frc$.
\end{theorem}
\begin{proof}
Because of item (i), $\frg_0=\espan{\bar h_1,\ldots,\bar h_n}$
holds, and condition (ii) implies that the dimension of $\frg_0$
equals the rank of $A$ (any linear combination of the rows of $A$
which gives $0$ induces a linear combination of the $\bar h_i$'s
which is central, and hence spans an ideal with trivial intersection
with $\frg_{-1}\oplus\frg_1$).

As in the proof of Theorem \ref{th:derivedgAtau} take a vector space
$\hat\frg_0$ of dimension $n-\rank A$, and elements $\bar
\alpha_1,\ldots,\bar\alpha_n\in (\frg_0\oplus\hat\frg_0)^*$ such
that
\[
\left\{\begin{aligned} &\bar\alpha_j(\bar h_i)=a_{ij}\ \text{for any
$i,j$,}\\
 &\text{$\bar\alpha_1,\ldots,\bar\alpha_n$ are linearly independent.}
 \end{aligned}\right.
\]

The Lie superalgebra $\bar \frg$ constructed as in the proof of
Theorem \ref{th:derivedgAtau} satisfies now the hypotheses of
Theorem \ref{th:centerlessgAtau}. Hence $\bar\frg$ is isomorphic to
the centerless contragredient Lie superalgebra $\frg(A,\tau)/\frc$,
and hence $\frg=[\bar\frg,\bar\frg]$ is isomorphic to the centerless
derived contragredient Lie superalgebra
$\frg'(A,\tau)/\frc=[\frg(A,\tau)/\frc,\frg(A,\tau)/\frc]$.
\end{proof}

\bigskip

In order to present the Lie superalgebras $\frg(S,S')$ in the
Extended Freudenthal's Magic Square as contragredient Lie
superalgebras, it will be useful first to describe them in terms
similar to those used in \cite{ElduqueMagicII}.

Let $V$ be, as before, a two dimensional vector space endowed with a
nonzero alternating bilinear form $\langle .\vert .\rangle$. For any
$n\in\bN$ consider copies $V_1,\ldots,V_n$ of $V$, and for any
subset $\sigma\subseteq \{1,\ldots,n\}$ take the module for
$\frsp(V)^n=\frsp(V_1)\oplus \cdots\oplus\frsp(V_n)$ given by:
\begin{equation}\label{eq:Vsigma}
V(\sigma)=\begin{cases}
 \frsp(V_1)\oplus\cdots\oplus \frsp(V_n)&\text{if
 $\sigma=\emptyset$,}\\
 V_{i_1}\otimes\cdots \otimes V_{i_r}&\text{if
 $\sigma=\{i_1,\ldots,i_r\}$, $1\leq i_1<\cdots<i_r\leq n$.}
 \end{cases}
\end{equation}
($V_i$ is the natural module for $\frsp(V_i)$ annihilated by
$\frsp(V_j)$ for $j\ne i$.)

Identify, as in \cite[Section 2]{ElduqueMagicII}, any subset
$\sigma\subseteq \{1,\ldots,n\}$ with the element
$(\sigma_1,\ldots,\sigma_n)\in\bZ_2^n$, such that $\sigma_i=1$ if
and only if $i\in\sigma$. Then for any $\sigma,\tau\in\bZ_2^n$,
consider the natural $\frsp(V)^n$-invariant map:
\begin{equation}\label{eq:phisigmatau}
\varphi_{\sigma,\tau}:V(\sigma)\times V(\tau)\longrightarrow
V(\sigma+\tau)
\end{equation}
defined as follows:

\begin{itemize}

\item
If $\sigma\ne \tau$ and $\sigma\ne\emptyset\ne\tau$, then
$\varphi_{\sigma,\tau}$ is obtained by contraction, by means of
$\langle .\vert .\rangle$ in the indices $i\in\sigma\cap\tau$
($\sigma_i=1=\tau_i$). Thus, for instance,
\[
\varphi_{\{1,2,3\},\{1,3,4\}}
 (u_1\otimes u_2\otimes u_3,v_1\otimes v_3\otimes v_4)=
 \langle u_1\vert v_1\rangle\langle u_3\vert v_3\rangle
     u_2\otimes v_4
\]
for any $u_1,v_1\in V_1$, $u_2\in V_2$, $u_3,v_3\in V_3$ and $v_4\in
V_4$.
\smallskip
\item
$\varphi_{\emptyset,\emptyset}$ is the Lie bracket in
$\frsp(V)^n=\frsp(V_1)\oplus\cdots\oplus\frsp(V_n)$.
\smallskip
\item
For any $\sigma\ne\emptyset$, $\varphi_{\emptyset,\sigma} =
-\varphi_{\sigma,\emptyset}$ is given by the natural action of the
Lie algebra $\frsp(V)^n$ on $V(\sigma)$. Thus, for instance,
\[
\varphi_{\emptyset,\{1,3\}}\bigl((s_1,\ldots,s_n),
 u_1\otimes u_3\bigr)=s_1(u_1)\otimes u_3 +
   u_1\otimes s_3(u_3),
\]
for any $s_i\in \frsp(V)$, $i=1,\ldots,n$, and $u_1\in V_1$, $u_3\in
V_3$.
\smallskip
\item
Finally, for any $\sigma=\{i_1,\ldots,i_r\}\ne \emptyset$,
$\varphi_{\sigma,\sigma}$ is given by:
\[
\varphi_{\sigma,\sigma} (u_{i_1}\otimes \cdots\otimes u_{i_r},
   v_{i_1}\otimes \cdots \otimes v_{i_r})=
   \sum_{j=1}^r\Bigl(
   \prod_{k\ne j}\langle u_{i_k}\vert v_{i_k}\rangle
   \nu_{i_j}(\gamma_{u_{i_j},v_{i_j}})\Bigr)\,,
\]
for any $u_{i_j},v_{i_j}\in V_{i_j}$, $j=1,\ldots,r$, where
$\nu_i:\frsp(V_i)\rightarrow \frsp(V)^n$ denotes the canonical
inclusion into the $i^{\text{\small th}}$-component, and
$\gamma_{u,v}$ has been defined in \eqref{eq:gammas}.
\end{itemize}

\medskip

\begin{example}\label{ex:d21}
The description of $\frd_{2,1}$ (characteristic $3$) in
\eqref{eq:d21} shows that, with $n=3$,
\begin{equation}\label{eq:d21Vs}
\frd_{2,1}=V(\emptyset)\oplus V(\{1,2,3\}),
\end{equation}
with even part $V(\emptyset)=\frsp(V_1)\oplus
\frsp(V_2)\oplus\frsp(V_3)$, and odd part $V(\{1,2,3\})=V_1\otimes
V_2\otimes V_3$, and where
\[
[x_\sigma,y_\tau]
=\epsilon_{\frd_{2,1}}(\sigma,\tau)\varphi_{\sigma,\tau}(x_\sigma,y_\tau),
\]
with
\begin{equation}\label{eq:epsilond21}
\begin{gathered}
\epsilon_{\frd_{2,1}}(\emptyset,\emptyset)
 =\epsilon_{\frd_{2,1}}(\emptyset,\{1,2,3\})
 =\epsilon_{\frd_{2,1}}(\{1,2,3\},\emptyset)=1,\\
 \epsilon_{\frd_{2,1}}(\{1,2,3\},\{1,2,3\})=-1.
 \end{gathered}
\end{equation}
(Same behavior as $\frd_4$ in \cite[(2.18) and
(2.19)]{ElduqueMagicII}.)

Take the canonical basic elements $\epsilon_1=(1,0,0)$,
$\epsilon_2=(0,1,0)$ and $\epsilon_3=(0,0,1)$ in $\bZ^3$, and for
any $i=1,2,3$ consider a symplectic basis $\{v_i,w_i\}$ of $V_i$
(that is, $\langle v_i\vert w_i\rangle =1$), and the basic elements
in each $\frsp(V_i)$:
\begin{equation}\label{eq:hieifi}
\begin{split}
h_i&=\gamma_{v_i,w_i}:\begin{cases} v_i\mapsto -v_i,\\
   w_i\mapsto w_i, \end{cases} \\
e_i&=\gamma_{w_i,w_i}:\begin{cases} v_i\mapsto -2w_i=w_i,\\
   w_i\mapsto 0, \end{cases} \\
f_i&=-\gamma_{v_i,v_i}:\begin{cases} v_i\mapsto 0,\\
   w_i\mapsto -2v_i=v_i, \end{cases}
\end{split}
\end{equation}
which satisfy
\[
[h_i,e_i]=2e_i,\quad [h_i,f_i]=-2f_i,\quad [e_i,f_i]=h_i.
\]

The Lie superalgebra $\frd_{2,1}$ is $\bZ^3$-graded by assigning
$\deg w_i=\epsilon_i=-\deg v_i$, $i=1,2,3$. Hence, $\deg h_i=0$,
$\deg e_i=2\epsilon_i$, $\deg f_i=-2\epsilon_i$, for any $i$.

Then,
\begin{equation}\label{eq:Phid21}
\Phi_{\frd_{2,1}}=\{\pm 2\epsilon_i: i=1,2,3\}\cup\{\pm
\epsilon_1\pm\epsilon_2\pm\epsilon_3\}
\end{equation}
is the set of the nonzero degrees that appear in $\frd_{2,1}$.
Moreover, the subalgebra $\frh=kh_1\oplus kh_2\oplus kh_3$ is a
Cartan subalgebra of $\frd_{2,1}$ and there is a natural
homomorphism of abelian groups:
\[
\begin{split}
R:\bZ^3&\longrightarrow \frh^*\\
\epsilon_i&\mapsto R(\epsilon_i)(:h_j\mapsto \delta_{ij}).
\end{split}
\]
The image under $R$ of $\Phi_{\frd_{2,1}}$ is precisely the set of
roots of $\frh$ in $\frd_{2,1}$.

Consider the lexicographic order on $\bZ^3$ with
$\epsilon_1>\epsilon_2>\epsilon_3>0$. The set of the positive
elements in $\Phi_{\frd_{2,1}}$ which are not sums of two positive
elements is
\begin{equation}\label{eq:Pid21}
\Pi_{\frd_{2,1}}=\{\alpha_1=2\epsilon_2,\alpha_2=\epsilon_1-\epsilon_2-\epsilon_3,
 \alpha_3=2\epsilon_3\},
\end{equation}
whose elements are linearly independent over $\bZ$ and satisfy that
$\Phi_{\frd_{2,1}}$ is contained in
$\bZ\Pi_{\frd_{2,1}}=\bZ\alpha_1\oplus\bZ\alpha_2\oplus\bZ\alpha_3$.

Let $Q$ be the free abelian group $\bZ\Pi_{\frd_{2,1}}$, take
$\tau=\{2\}$ and consider the following elements in $\frd_{2,1}$:
\[
\begin{aligned}
E_1&=e_2, & E_2&=w_1\otimes v_2\otimes v_3, & E_3&=e_3,\\
F_1&=f_2, & F_2&=-v_1\otimes w_2\otimes w_3, & F_3&=f_3,\\
H_1&=h_2, & H_2&=h_1-h_2-h_3, & H_3&=h_3.
\end{aligned}
\]
Then $\bigl(\frh,\{
R(\alpha_1),R(\alpha_2),R(\alpha_3)\},\{H_1,H_2,H_3\}\bigr)$ is a
realization of the matrix:
\begin{equation}\label{eq:Ad21}
A_{\frd_{2,1}}=\begin{pmatrix} 2&-1&0\\ 1&0&1\\ 0&-1&2\end{pmatrix},
\end{equation}
which is the Cartan matrix $D_1$ in \cite[p.~55]{KacLieSuper}, and
which corresponds to the Lie superalgebra $D(2,1)\cong D(2,1;1)$,
that is, to the orthosymplectic Lie superalgebra $\frosp(4,2)$.
\end{example}

The Lie superalgebra $\frd_{2,1}$ is $\bZ$-graded with
$(\frd_{2,1})_0=\frh$, $(\frd_{2,1})_1=kE_1\oplus kE_2\oplus kE_3$,
and $(\frd_{2,1})_{-1}=kF_1\oplus kF_2\oplus kF_3$. It is easily
checked to be simple. Now, Theorem \ref{th:gAtau}, with the elements
$H_i,E_i,F_i$ above, $i=1,2,3$, gives:

\begin{proposition}\label{pr:d21}
The Lie superalgebra $\frd_{2,1}$ is isomorphic to the
contragredient Lie superalgebra
$\frg\bigl(A_{\frd_{2,1}},\{2\}\bigr)$.
\end{proposition}

\bigskip

\section{The Lie superalgebras in the Extended Freudenthal's Magic
Square}

In \cite{ElduqueMagicII}, most of the Lie algebras in Freudenthal's
Magic Square were described, in the split case, by means of the
$V(\sigma)$'s and $\varphi_{\sigma,\tau}$'s of the previous section.
This is possible too in the superalgebra setting.

This section will be devoted to get such descriptions, or a similar
one, for each Lie superalgebra in the Extended Freudenthal's Magic
Square, constructed from a couple of split symmetric composition
(super)algebras. This will be used to find a description of all
these superalgebras as contragredient Lie superalgebras.

Throughout this section, the characteristic of the ground field $k$
will always be assumed to be $3$.

\bigskip

\subsection{$\frg(S_1,S_{4,2})$}\label{sub:gS1S42}

The symmetric composition superalgebra $(S_1,\bullet,q)$ is just
$k1$, with $1\bullet 1=1$ and $q(1)=1$ (so $b(1,1)=2$). Thus
$\tri(S_1)=0$ and \eqref{eq:gSSprime}, \eqref{eq:d21Vs}  and
\eqref{eq:iotaS42} show that
\[
\begin{split}
\frg(S_1,S_{4,2})&=\tri(S_{4,2})\oplus \iota_0(S_1\otimes
S_{4,2})\oplus
\iota_1(S_1\otimes S_{4,2})\oplus\iota_2(S_1\otimes S_{4,2})\\
 &=\bigl(V(\emptyset)\oplus
 V(\{1,2,3\})\bigr)\oplus\bigl(V(\{1,2\})\oplus V(\{3\})\bigr)\\
 &\null\qquad \oplus\bigl(V(\{2,3\})\oplus V(\{1\})\bigr)
 \oplus\bigl(V(\{1,3\})\oplus V(\{2\})\bigr)\\
 &=\bigoplus_{\sigma\in 2^{\{1,2,3\}}}V(\sigma),
\end{split}
\]
where $2^{\{1,2,3\}}$ denotes the power set of $\{1,2,3\}$. Besides,
\[
\frg(S_1,S_{4,2})\subo=\bigoplus_{\substack{\sigma\in
2^{\{1,2,3\}}\\ \lvert \sigma\rvert\ \text{even}}} V(\sigma)\qquad
\text{and}\qquad\frg(S_1,S_{4,2})\subuno=\bigoplus_{\substack{\sigma\in
2^{\{1,2,3\}}\\ \lvert \sigma\rvert\ \text{odd}}} V(\sigma).
\]
By invariance of the Lie bracket under the action of
$\frsp(V_1)\oplus\frsp(V_2)\oplus\frsp(V_3)=V(\emptyset)$, it
follows that
\begin{equation}\label{eq:epsilon}
[x_\sigma,y_\tau]=\epsilon(\sigma,\tau)\varphi_{\sigma,\tau}(x_\sigma,y_\tau),
\end{equation}
for any $\sigma,\tau\in 2^{\{1,2,3\}}$, for a suitable map
\[
\epsilon: 2^{\{1,2,3\}}\times 2^{\{1,2,3,\}}\rightarrow k^\times.
\]
The multiplication $\iota_i(1\otimes S_{4,2})\times
\iota_{i+1}(1\otimes S_{4,2})\rightarrow \iota_{i+2}(1\otimes
S_{4,2})$ is given in \eqref{eq:bulletiotaS42}, and for any $x,y\in
S_{4,2}$ (see \eqref{eq:txy}):
\[
[\iota_0(1\otimes x),\iota_0(1\otimes y)]=2t_{x,y}
=\bigl(2\sigma_{x,y},b(x,y)1-2r_xl_y,b(x,y)1-2l_xr_y\bigr),
\]
as $b(1,1)=2$ in $S_1$. This element is identified with the element
$2\rho^{-1}(\sigma_{x,y})$ in \eqref{eq:rhod21}, which is given by
\eqref{eq:sigmasS42}. By cyclic symmetry, one completes the
information about the map $\epsilon$ in \eqref{eq:epsilon},
displayed in Table \ref{ta:epsilonS42}, which is exactly Table 2 in
\cite{ElduqueMagicII}, corresponding to $\frg(S_1,S_8)$, but with
the index $4$ taken out. From this description it readily follows
that $\frg(S_1,S_{4,2})$ is simple.

\begin{table}[h!]
$$
\vbox{\offinterlineskip \halign{\hfil$#$\hfil\enspace\vreglon
 &\enspace\hfil$#$\hfil\enspace\vregla
 &\enspace\hfil$#$\hfil\enspace\vregla
 &\enspace\hfil$#$\hfil\enspace\vregla
 &\enspace\hfil$#$\hfil\enspace&\hskip .2pt\vregla#
 &&\enspace\hfil$#$\hfil\enspace\vregla\cr
 \omit\enspace\hfil\raise2pt\hbox{$\epsilon$}%
  \enspace\hfil\vrule width1pt depth 4pt
   &\enspace\emptyset\enspace&\{1,2\}&\{2,3\}&\{1,3\}&\omit\vregleta
     &\omit$\{1,2,3\}$\vregla&\{3\}&\{1\}&\{2\}\cr
 \noalign{\hreglon}
 \emptyset&1&1&1&1&\omit\vregleta&1&1&1&1\cr
 \noalign{\hregla}
 \{1,2\}&1&-2&1&1&\omit\vregleta&1&-2&-1&-1\cr\noalign{\hregla}
 \{2,3\}&1&1&-2&1&\omit\vregleta&1&-1&-2&-1\cr\noalign{\hregla}
 \{1,3\}&1&1&1&-2&\omit\vregleta&1&-1&-1&-2\cr
 \multispan6{\hregletafill}&\multispan4{\hreglafill}\cr
 \{1,2,3\}&1&-1&-1&-1&&-1&1&1&1\cr\noalign{\hregla}
 \{3\}&1&2&-1&-1&&-1&-2&-1&-1\cr\noalign{\hregla}
 \{1\}&1&-1&2&-1&&-1&-1&-2&-1\cr\noalign{\hregla}
 \{2\}&1&-1&-1&2&&-1&-1&-1&-2\cr
 \noalign{\hregla}
}}
$$
\medskip
\caption{$\epsilon_{S_1,S_{4,2}}$}\label{ta:epsilonS42}
\end{table}

\medskip

In general, in order to obtain a description of the Lie superalgebra
$\frg(S,S')$ in the Extended Freudenthal's Magic Square as
contragredient Lie superalgebras, the complete description of the
map $\epsilon$ in \eqref{eq:epsilon} is not necessary and will not
be given. A detailed description of these maps appear in
\cite{Cu06}.

\smallskip

Now, with the notations and conventions in Example \ref{ex:d21}, the
Lie superalgebra $\frg=\frg(S_1,S_{4,2})$ is $\bZ^3$-graded, and the
set of nonzero degrees that appear is
\[
\Phi_{S_1,S_{4,2}}=\bigl\{\pm
2\epsilon_i,\pm\epsilon_i\pm\epsilon_j,
\pm\epsilon_1\pm\epsilon_2\pm\epsilon_3,\pm\epsilon_i: 1\leq i<j\leq
3\bigr\}.
\]
The subalgebra $\frh=kh_1\oplus kh_2\oplus kh_3$ (recall
\eqref{eq:hieifi} that $h_i=\gamma_{v_i,w_i}$ belongs to
$\frsp(V_i)\subseteq V(\emptyset)$, where $\{v_i,w_i\}$ is a fixed
symplectic basis of $V_i$) is a Cartan subalgebra of $\frg$ and the
image of $\Phi_{S_1,S_{4,2}}$ under the natural homomorphism of
abelian groups
\[
\begin{split}
R:\bZ^3&\longrightarrow \frh^*\\
 \epsilon_i&\mapsto R(\epsilon_i)(:h_j\mapsto \delta_{ij})
\end{split}
\]
is precisely the set of roots of $\frg$ relative to $\frh$.

Consider the lexicographic order on $\bZ^3$ with
$\epsilon_1>\epsilon_2>\epsilon_3>0$. The set of positive elements
in $\Phi_{S_1,S_{4,2}}$ which are not sums of positive elements (the
\emph{irreducible} degrees) is
\[
\Pi=\{\alpha_1=\epsilon_2-\epsilon_3,\alpha_2=\epsilon_3,
 \alpha_3=\epsilon_1-\epsilon_2-\epsilon_3\},
\]
which is a $\bZ$-basis of $\bZ^3$.

Then $\bigl(\frh,\{R(\alpha_1),R(\alpha_2),R(\alpha_3)\},
\{h_2-h_3,2h_3,-h_1+h_2+h_3\}\bigr)$ is a realization of the matrix
\[
A_{S_1,S_{4,2}}=\begin{pmatrix} 2&-1&0\\ 1&2&1\\ 0&1&0\end{pmatrix},
\]
with associated Dynkin diagram (using the conventions in
\cite[Tables IV and V]{KacLieSuper})
\[
\SunoScuatrodos
\]

\smallskip

With $\tau=\{2,3\}$, consider the generators
\[
\begin{aligned}
H_1&=h_2-h_3, & H_2&=2h_3, & H_3&=-h_1+h_2+h_3,\\
E_1&=w_2\otimes v_3, & E_2&=w_3, & E_3&=w_1\otimes v_2\otimes v_3,\\
F_1&=v_2\otimes w_3, & F_2&=-v_3, & F_3&=-v_1\otimes w_2\otimes w_3.
\end{aligned}
\]
As for Proposition \ref{pr:d21}, Theorem \ref{th:gAtau} implies the
following description of $\frg(S_1,S_{4,2})$:

\begin{proposition}\label{pr:gS1S42}
The Lie superalgebra $\frg(S_1,S_{4,2})$ is isomorphic to the
contragredient Lie superalgebra
$\frg\bigl(A_{S_1,S_{4,2}},\{2,3\}\bigr)$.
\end{proposition}

\smallskip

Moreover, the set of even and odd nonzero degrees are
\[
\begin{split}
\bigl(\Phi_{S_1,S_{4,2}}\bigr)\subo&=
 \{\pm 2\epsilon_i, \pm\epsilon_i\pm\epsilon_j: 1\leq i<j\leq 3\},\\
\bigl(\Phi_{S_1,S_{4,2}}\bigr)\subuno&=
 \{\pm \epsilon_i, \pm\epsilon_1\pm\epsilon_2\pm\epsilon_3: 1\leq i\leq
 3\}.
\end{split}
\]
The set of ``irreducible''even degrees for the lexicographic order
considered is
\[
\Pi\subo=\{\beta_1=\epsilon_1-\epsilon_2,\beta_2=\epsilon_2-\epsilon_3,
 \beta_3=2\epsilon_3\},
\]
and it can be concluded from here that the even part
$\frg\subo=\frg(S_1,S_{4,2})\subo$ is isomorphic to the
contragredient Lie algebra with Cartan matrix
\[
\begin{pmatrix} 2&-1&0\\ -1&2&-2\\ 0&-1&2\end{pmatrix}.
\]
That is, $\frg(S_1,S_{4,2})\subo$ is isomorphic to the symplectic
Lie algebra $\frsp_6$ (type $C_3$). Here the appropriate basis of
$\frh$ is
\[
\{\tilde H_1=h_1-h_2,\tilde H_2=h_2-h_3,\tilde H_3=h_3\}.
\]
Also, the set of simple roots $R(\Pi\subo)$ of $\frg\subo$ induces a
triangular decomposition $\frg\subo=\frn^-\oplus\frh\oplus\frn^+$,
where
$\frn^+=\oplus_{0<\alpha\in(\Phi_{S_1,S_{4,2}})\subo}(\frg\subo)_{R(\alpha)}$
is the sum of root spaces (nonzero homogeneous components in the
$\bZ^3$-grading) corresponding to the positive degrees in
$(\Phi_{S_1,S_{4,2}})\subo$, and similarly for $\frn^-$. The
$\frg\subo$-module $\frg\subuno$ is $\bZ^3$-graded consistently with
the action of $\frg\subo$, it is easily seen to be irreducible, and
to contain the highest weight vector $w_1\otimes w_2\otimes w_3$,
that is, $[\frn^+,w_1\otimes w_2\otimes w_3]=0$. Its highest weight
is $R(\epsilon_1+\epsilon_2+\epsilon_3)=\omega_3$, which satisfies
$\omega_3(\tilde H_1)=0=\omega_3(\tilde H_2)$, $\omega_3(\tilde
H_3)=1$.

Since the proof of the uniqueness result in \cite[Theorem A in \S
20.3]{Humphreys} remains valid in this setting (note that instead of
grading over the $\bZ$-linear combinations of weights, which are
elements of $\frh^*$, we grade over a true lattice $\bZ^3$),
$\frg\subuno$ is the unique $\frg\subo$-module with such a highest
weight. Denote it by $V(\omega_3)$. Note that $\dim\frg\subuno=14$.

Let us give the most natural presentation of this module. Consider
the matrix Lie algebra $\frsp_6$, and its natural six-dimensional
module, which is endowed with an alternating invariant bilinear form
$\{.\vert.\}$. Let $\{a_1,a_2,a_3,b_1,b_2,$ $b_3\}$ be a symplectic
basis (that is, $\{a_i\vert b_j\}=\delta_{ij}$, $\{a_i\vert
a_j\}=0=\{b_i\vert b_j\}$ for any $i,j$). With
$\gamma_{x,y}=\{x\vert .\}y+\{y\vert .\}x$, the subspace spanned by
$\gamma_{a_i,b_i}$, $i=1,2,3$ is a Cartan subalgebra of
$\frsp(W)\simeq\frsp_6$ and the weights of $W$ are $\{\pm\delta_i:
i=1,2,3\}$, where $\delta_i(\gamma_{a_j,b_j})=\delta_{ij}$ for any
$i,j$.

Consider the $\frsp_6$-module $\bigwedge^3 W$ and the homomorphism
\[
\begin{split}
\varphi: \bigwedge^3W&\longrightarrow W\\
z_1\wedge z_2\wedge z_3&\mapsto \{z_1\vert z_2\}z_3+\{z_2\vert
z_3\}z_1+\{z_3\vert z_1\}z_2.
\end{split}
\]
Then $\dim\ker\varphi=14$, the weights of $\ker\varphi$ are
$\pm\delta_1\pm\delta_2\pm\delta_3$ and $\pm\delta_i$ ($i=1,2,3$),
all of them of multiplicity $1$. The element $b_1\wedge b_2\wedge
b_3$ is a highest weight vector of weight
$\delta_1+\delta_2+\delta_3=\omega_3$ and $\ker\varphi$ is
irreducible. Therefore, up to isomorphism $V(\omega_3)=\ker\varphi$.

Let us summarize the above discussion:

\begin{proposition}\label{pr:gS1S42evenodd}
The Lie superalgebra $\frg(S_1,S_{4,2})$ is simple with even part
isomorphic to the symplectic Lie algebra $\frsp_6$ and odd part
isomorphic to the irreducible module of dimension $14$ above.
\end{proposition}

There is no counterpart in characteristic $0$ (see
\cite{KacLieSuper}) to this simple Lie superalgebra.

\bigskip

\subsection{$\frg(S_4,S_{4,2})$}

The symmetric composition superalgebra $S_4$ is the even part of
$S_{4,2}$, and its triality Lie algebra is given too by the even
part of the triality Lie superalgebra of $S_{4,2}$. Therefore, six
copies of $V$ are needed: $V_1$, $V_2$ and $V_3$ for $S_4$, and
$V_4$, $V_5$ and $V_6$ for $S_{4,2}$.

With the same sort of arguments used so far, one has:
\[
\begin{split}
\frg(S_4,S_{4,2})
 &=\bigl(\tri(S_4)\oplus\tri(S_{4,2})\bigr)\oplus
  \bigl(\oplus_{i=0}^2\iota_i(S_4\otimes S_{4,2})\bigr)\\
 &=\oplus_{\sigma\in\calS_{S_4,S_{4,2}}}V(\sigma),
\end{split}
\]
where
\begin{multline*}
\calS=\calS_{S_4,S_{4,2}}=\bigl\{
\emptyset,\{4,5,6\},\\
 \{1,2,4,5\},\{1,2,6\},\{2,3,5,6\},\{2,3,4\},\{1,3,4,6\},\{1,3,5\}\bigr\}.
\end{multline*}
For instance,
\[
\iota_1(S_4\otimes S_{4,2})=(V_2\otimes V_3)\otimes\bigl((V_5\otimes
V_6)\oplus V_4\bigr)=V\bigl(\{2,3,5,6\}\bigr)\oplus
V\bigl(\{2,3,4\}\bigr).
\]

Like in all the other cases, the even (respectively odd) part is the
sum of the $V(\sigma)$'s with $\sigma$  containing an even (resp.
odd) number of elements.

The multiplication presents the form in \eqref{eq:epsilon} for a
suitable map $\epsilon:\calS\times\calS\rightarrow k^\times$. Here
the nonzero even and odd degrees are:
\[
\begin{split}
\Phi\subo&=\{\pm2\epsilon_i: 1\leq i\leq 6\}\\
 &\null\qquad \cup
 \{\pm\epsilon_1\pm\epsilon_2\pm\epsilon_4\pm\epsilon_5,
 \pm\epsilon_2\pm\epsilon_3\pm\epsilon_5\pm\epsilon_6,
 \pm\epsilon_1\pm\epsilon_3\pm\epsilon_4\pm\epsilon_6\},\\[4pt]
\Phi\subuno&=\{\pm\epsilon_4\pm\epsilon_5\pm\epsilon_6,
   \pm\epsilon_1\pm\epsilon_2\pm\epsilon_6,
   \pm\epsilon_2\pm\epsilon_3\pm\epsilon_4,
   \pm\epsilon_1\pm\epsilon_3\pm\epsilon_5\}.
\end{split}
\]
With the lexicographic order given by $0<\epsilon_1<\cdots
<\epsilon_6$, the set of irreducible degrees is
\begin{multline*}
\Pi=\{\alpha_1=2\epsilon_1,\alpha_2=\epsilon_5-\epsilon_4-\epsilon_2-\epsilon_1,
\alpha_3=2\epsilon_2,\\
 \alpha_4=\epsilon_4-\epsilon_3-\epsilon_2,\alpha_5=2\epsilon_3,
 \alpha_6=\epsilon_6-\epsilon_5-\epsilon_4\},
\end{multline*}
which is a linearly independent set over $\bZ$. Then $\Phi$ is
contained in $Q=\bZ\Pi$ (which is isomorphic to $\bZ^6$), so that
any positive element in $\Phi$ is a sum of elements in $\Pi$.
Consider the elements (same conventions as before)
\begin{multline*}
H_1=h_1,\ H_2=-h_5+h_4+h_2+h_1,\ H_3=h_2,\ H_4=h_4-h_3-h_2,\\
H_5=h_3,\ H_6=h_6-h_5-h_4,
\end{multline*}
which give a realization
\[
\bigl(\frh,\{R(\alpha_1),R(\alpha_2),R(\alpha_3),R(\alpha_4),R(\alpha_5),R(\alpha_6)\},
\{H_1,H_2,H_3,H_4,H_5,H_6\}\bigr)
\]
of the regular matrix (with associated Dynkin diagram):
\[
A_{S_4,S_{4,2}}=\begin{pmatrix} 2&-1&0&0&0&0\\ -1&2&-1&0&0&0\\
 0&-1&2&-1&0&0\\ 0&0&1&0&1&-1 \\ 0&0&0&-1&2&0 \\
 0&0&0&-1&0&0\end{pmatrix},\qquad \ScuatroScuatrodos
\]
Then:

\begin{proposition}\label{pr:gS4S42}
The Lie superalgebra $\frg(S_4,S_{4,2})$ is isomorphic to the
contragredient Lie superalgebra
$\frg\bigl(A_{S_4,S_{4,2}},\{4,6\}\bigr)$.
\end{proposition}

Also, the set of irreducible even degrees is:
\begin{multline*}
\Pi\subo=\{\beta_1=2\epsilon_3,\beta_2=\epsilon_6-\epsilon_5-\epsilon_3-\epsilon_2,
\beta_3=2\epsilon_2,\\
 \beta_4=\epsilon_5-\epsilon_4-\epsilon_2-\epsilon_1,\beta_5=2\epsilon_1,
 \beta_6=2\epsilon_4\},
\end{multline*}
with associated Cartan matrix of type $D_6$. The odd part is an
irreducible module for the even part, with highest weight vector
$R(\epsilon_4+\epsilon_5+\epsilon_6)=\omega_6$, and the same sort of
arguments used for $\frg(S_1,S_{4,2})$ gives:

\begin{proposition}\label{pr:gS4S42evenodd}
The Lie superalgebra $\frg(S_4,S_{4,2})$ is simple with even part
isomorphic to the orthogonal Lie algebra $\frso_{12}$, and odd part
isomorphic to the spin module for $\frso_{12}$.
\end{proposition}

Note that the even part is just $\frg(S_4,S_4)$, which was
determined in \cite{ElduqueMagicII}.

\begin{corollary}
The Lie superalgebra $\frg(S_4,S_{4,2})$ is isomorphic to the Lie
superalgebra in \cite[theorem 3.2(v)]{ElduqueNewSimple3} and
\cite[Theorem 4.1(ii)\, $(l=6)$]{ElduqueNewSimple35}.
\end{corollary}

\bigskip

\subsection{$\frg(S_8,S_{4,2})$}

Here, four copies of $V$ are needed for $S_8$ (see
\cite[(2.14)]{ElduqueMagicII}) and three copies for $S_{4,2}$. The
indices $1$, $2$, $3$ and $4$ will be used for $S_8$, while $5$, $6$
and $7$ will be reserved for $S_{4,2}$. Then
\[
\begin{split}
\frg(S_8,S_{4,2})
 &=\bigl(\tri(S_8)\oplus\tri(S_{4,2})\bigr)\oplus
  \bigl(\oplus_{i=0}^2\iota_i(S_8\otimes S_{4,2})\bigr)\\
 &=\oplus_{\sigma\in\calS_{S_8,S_{4,2}}}V(\sigma),
\end{split}
\]
where
\[
\begin{split}
\calS=\calS_{S_8,S_{4,2}}=\bigl\{
\emptyset,&\{1,2,3,4\},\{5,6,7\},\\
&\quad \{1,2,5,6\},\{3,4,5,6\},\{1,2,7\},\{3,4,7\},\\
&\quad \{2,3,6,7\},\{1,4,6,7\},\{2,3,5\},\{1,4,5\},\\
&\quad \{1,3,5,7\},\{2,4,5,7\},\{1,3,6\},\{2,4,6\}\bigr\}.
\end{split}
\]
Hence,
\[
\begin{split}
\Phi\subo&=\bigl\{\pm 2\epsilon_i: i=1,\ldots,7\bigr\}\cup\bigl\{
  \pm\epsilon_1\pm\epsilon_2\pm\epsilon_3\pm\epsilon_4,\\
  &\qquad \pm\epsilon_1\pm\epsilon_2\pm\epsilon_5\pm\epsilon_6,
   \pm\epsilon_3\pm\epsilon_4\pm\epsilon_5\pm\epsilon_6,
   \pm\epsilon_2\pm\epsilon_3\pm\epsilon_6\pm\epsilon_7,\\
  &\qquad \pm\epsilon_1\pm\epsilon_4\pm\epsilon_6\pm\epsilon_7,
   \pm\epsilon_1\pm\epsilon_3\pm\epsilon_5\pm\epsilon_7,
   \pm\epsilon_2\pm\epsilon_4\pm\epsilon_5\pm\epsilon_7\bigr\},\\[6pt]
\Phi\subuno&=\bigl\{\pm\epsilon_5\pm\epsilon_6\pm\epsilon_7,
     \pm\epsilon_1\pm\epsilon_2\pm\epsilon_7,
     \pm\epsilon_3\pm\epsilon_4\pm\epsilon_7,
     \pm\epsilon_2\pm\epsilon_3\pm\epsilon_5,\\
   &\qquad \pm\epsilon_1\pm\epsilon_4\pm\epsilon_5,
      \pm\epsilon_1\pm\epsilon_3\pm\epsilon_6,
      \pm\epsilon_2\pm\epsilon_4\pm\epsilon_6\bigr\},
\end{split}
\]
and with the lexicographic order with $0<\epsilon_1<\cdots
<\epsilon_7$:
\begin{multline*}
\Pi=\bigl\{\alpha_1=\epsilon_6-\epsilon_5-\epsilon_4-\epsilon_3,
   \alpha_2=2\epsilon_2,\alpha_3=2\epsilon_3,\\
   \alpha_4=\epsilon_4-\epsilon_3-\epsilon_2-\epsilon_1,
    \alpha_5=2\epsilon_1,\alpha_6=\epsilon_5-\epsilon_4-\epsilon_1,
   \alpha_7=\epsilon_7-\epsilon_6-\epsilon_5\bigr\},
\end{multline*}
which is a $\bZ$-linearly independent set with $\Phi\subseteq
\bZ\Pi$, so that any positive element in $\Phi$ is a sum of elements
in $\Pi$. The associated matrix and Dynkin diagram are here:
\[
A_{S_8,S_{4,2}}=\begin{pmatrix} 2&0&-1&0&0&0&0\\ 0&2&0&-1&0&0&0\\
  -1&0&2&-1&0&0&0\\ 0&-1&-1&2&-1&0&0\\ 0&0&0&-1&2&-1&0\\
  0&0&0&0&1&0&-1\\ 0&0&0&0&0&-1&0 \end{pmatrix},\qquad
\SochoScuatrodos
\]

\smallskip

\begin{proposition}\label{pr:gS8S42}
The Lie superalgebra $\frg(S_8,S_{4,2})$ is isomorphic to the
contragredient Lie superalgebra
$\frg\bigl(A_{S_8,S_{4,2}},\{6,7\}\bigr)$.
\end{proposition}

Also, in this situation:
\[
\begin{split}
\Pi\subo&=\bigl\{
\beta_1=\epsilon_7-\epsilon_6-\epsilon_4-\epsilon_1,
\beta_2=2\epsilon_2, \beta_3=2\epsilon_1,\\
 &\qquad\beta_4=\epsilon_4-\epsilon_3-\epsilon_2-\epsilon_1,
 \beta_5=2\epsilon_3,\beta_6=\epsilon_6-\epsilon_5-\epsilon_4-\epsilon_3,
 \beta_7=2\epsilon_5\bigr\},
\end{split}
\]
whose associated Cartan matrix is of type $E_7$ (this also follows
since $\frg(S_8,S_{4,2})\subo$ is isomorphic to $\frg(S_8,S_4)$,
which was computed in \cite{ElduqueMagicII}).

The odd part, which has dimension $7\times 8=56$, is an irreducible
module with highest weight
$R(\epsilon_7+\epsilon_6+\epsilon_5)=\omega_7$. Thus,

\begin{proposition}\label{pr:gS8S42evenodd}
$\frg(S_8,S_{4,2})$ is a simple Lie superalgebra, whose even part is
isomorphic to the split simple Lie algebra of type $E_7$, and whose
odd part is its $56$-dimensional irreducible module $V(\omega_7)$.
\end{proposition}

\bigskip

\subsection{$\frg(S_{4,2},S_{4,2})$}

For $\frg(S_{4,2},S_{4,2})$, the indices $1$, $2$ and $3$ will refer
to the copies of $V$ associated to the first copy of $S_{4,2}$,
while the indices $4$, $5$ and $6$ will refer to the copies of $V$
related to the second copy of $S_{4,2}$. Then:
\[
\begin{split}
\frg(S_{4,2},S_{4,2})&=\bigl(\tri(S_{4,2})\oplus\tri(S_{4,2})\bigr)\oplus
 \bigl(\oplus_{i=0}^2\iota_i(S_{4,2}\otimes S_{4,2})\bigr)\\
 &=\oplus_{\sigma\in\calS_{S_{4,2},S_{4,2}}} V(\sigma),
\end{split}
\]
where
\[
\begin{split}
\calS_{S_{4,2},S_{4,2}}=\bigl\{\emptyset,&\{1,2,3\},\{4,5,6\},\\
  &\{1,2,4,5\},\{ 3,6\},\{1,2,6\},\{3,4,5\},\\
  &\{2,3,5,6,\},\{1,4\},\{2,3,4\},\{1,5,6\},\\
  &\{1,3,4,6\},\{2,5\},\{1,3,5\},\{2,4,6\}\bigr\}.
\end{split}
\]
Hence,
\[
\begin{split}
\Phi\subo&=\bigl\{\pm 2\epsilon_i: 1\leq i\leq 6\bigr\}\cup
 \bigl\{\pm\epsilon_1\pm\epsilon_2\pm\epsilon_4\pm\epsilon_5,
    \pm\epsilon_2\pm\epsilon_3\pm\epsilon_5\pm\epsilon_6,\\
    &\qquad \pm\epsilon_1\pm\epsilon_3\pm\epsilon_4\pm\epsilon_6,
    \pm\epsilon_3\pm\epsilon_6,\pm\epsilon_1\pm\epsilon_4,
    \pm\epsilon_2\pm\epsilon_5\bigr\},\\[6pt]
\Phi\subuno&=\bigl\{ \pm\epsilon_1\pm\epsilon_2\pm\epsilon_3,
    \pm\epsilon_4\pm\epsilon_5\pm\epsilon_6,
    \pm\epsilon_1\pm\epsilon_2\pm\epsilon_6,
    \pm\epsilon_3\pm\epsilon_4\pm\epsilon_5,\\
    &\qquad \pm\epsilon_2\pm\epsilon_3\pm\epsilon_4,
    \pm\epsilon_1\pm\epsilon_5\pm\epsilon_6,
    \pm\epsilon_1\pm\epsilon_3\pm\epsilon_5,
    \pm\epsilon_2\pm\epsilon_4\pm\epsilon_6\bigr\}.
\end{split}
\]
And with the lexicographic order with
$\epsilon_1>\cdots>\epsilon_6>0$,
\[
\begin{split}
\Pi&=\bigl\{\alpha_1=\epsilon_1-\epsilon_2-\epsilon_3,
  \alpha_2=\epsilon_3-\epsilon_4-\epsilon_5,\alpha_3=2\epsilon_5,\\
  &\qquad \alpha_4=\epsilon_4-\epsilon_5-\epsilon_6,
  \alpha_5=2\epsilon_6,
  \alpha_6=\epsilon_2-\epsilon_3-\epsilon_4\bigr\},
\end{split}
\]
which is a linearly independent set with $\Phi\subseteq \bZ\Pi$, so
that any positive element in $\Phi$ is a sum of elements in $\Pi$.
The associated matrix and Dynkin diagram are:
\[
A_{S_{4,2},S_{4,2}}=\begin{pmatrix}
 0&1&0&0&0&0\\ 1&0&-1&0&0&0\\ 0&-1&2&-1&0&0\\ 0&0&-1&0&-1&-1\\
  0&0&0&-1&2&0\\ 0&0&0&1&0&0\end{pmatrix},\qquad
  \ScuatrodosScuatrodos
\]

\begin{proposition}\label{pr:gS42S42}
The Lie superalgebra $\frg(S_{4,2},S_{4,2})$ is isomorphic to the
contragredient Lie superalgebra
$\frg\bigl(A_{S_{4,2},S_{4,2}},\{1,2,4,6\}\bigr)$.
\end{proposition}

Also,
\begin{multline*}
\Pi\subo=\bigl\{\beta_1=2\epsilon_4,
   \beta_2=\epsilon_1-\epsilon_2-\epsilon_4-\epsilon_5,
   \beta_3=2\epsilon_5,\\
    \beta_4=\epsilon_2-\epsilon_3-\epsilon_5-\epsilon_6,
   \beta_5=2\epsilon_6,
   \beta_6=\epsilon_3-\epsilon_6\bigr\},
\end{multline*}
with associated Cartan matrix of type $B_6$. The odd part is an
irreducible module for the even part with highest weight
$R(\epsilon_1+\epsilon_2+\epsilon_3)=\omega_6$, so by uniqueness, it
is the spin module for the even part:

\begin{proposition}\label{pr:S42S42evenodd}
$\frg(S_{4,2},S_{4,2})$ is a simple Lie superalgebra whose even part
is isomorphic to the orthogonal Lie algebra $\frso_{13}$ and whose
odd part is the spin module for its even part.
\end{proposition}

\begin{corollary}\label{co:S42S42}
The Lie superalgebra $\frg(S_{4,2},S_{4,2})$ is isomorphic to the
Lie superalgebra in \cite[Theorem 3.1(ii)\,
$(l=6)$]{ElduqueNewSimple35}.
\end{corollary}

\bigskip

\subsection{$\frg(S_{1,2},S_{4,2})$}

Here there is just one copy of $V$ involved in $S_{1,2}$, which will
carry index $1$, and three copies, with indices $2$, $3$ and $4$, in
$S_{4,2}$. Also, the triality Lie superalgebra $\tri(S_{1,2})$ will
be identified to the superalgebra $\frb_{0,1}$ in \eqref{eq:b01}.
Then
\[
\begin{split}
\frg(S_{1,2},S_{4,2})&=\bigl(\tri(S_{1,2})\oplus\tri(S_{4,2})\bigr)\oplus
  \bigl(\oplus_{i=0}^2\iota_i(S_{1,2}\otimes S_{4,2})\bigr)\\
  &=\oplus_{\sigma\in\calS_{S_{1,2},S_{4,2}}}V(\sigma),
\end{split}
\]
with
\[
\begin{split}
\calS_{S_{1,2},S_{4,2}}=\bigl\{\emptyset,&\{1\},\{2,3,4\},\{1,4\},\{1,2,3\},\{4\},\\
 &\{3,4\},\{1,2\},\{1,3,4\},\{2\},\{2,4\},\{1,3\},\{1,2,4\},\{3\}\bigr\},
\end{split}
\]
and
\[
\begin{split}
\Phi\subo&=\{\pm 2\epsilon_i:1\leq i\leq
   4\}\cup\{\pm\epsilon_i\pm\epsilon_j: 1\leq i<j\leq 4\},\\
\Phi\subuno&=\{\pm\epsilon_i\pm\epsilon_j\pm\epsilon_k: 1\leq
i<j<k\leq 4\}\cup\{\pm\epsilon_i:1\leq i\leq 4\}.
\end{split}
\]
With the lexicographic order
$\epsilon_1>\epsilon_2>\epsilon_3>\epsilon_4>0$,
\[
\Pi=\{\alpha_1=\epsilon_1-\epsilon_2-\epsilon_3,
 \alpha_2=\epsilon_3-\epsilon_4,
 \alpha_3=\epsilon_4,
 \alpha_4=\epsilon_2-\epsilon_3-\epsilon_4\},
\]
which is a linearly independent set with $\Phi\subseteq \bZ\Pi$, so
that any positive element in $\Phi$ is a sum of elements in $\Pi$.
The associated matrix and Dynkin diagram are:
\[
A_{S_{1,2},S_{4,2}}=\begin{pmatrix} 0&1&0&0\\ -1&2&-1&0\\
   0&-2&2&-2\\ 0&0&1&0\end{pmatrix},\qquad \SunodosScuatrodos
\]
and, therefore:

\begin{proposition}\label{pr:S12S42}
The Lie superalgebra $\frg(S_{1,2},S_{4,2})$ is isomorphic to the
contragredient Lie superalgebra
$\frg\bigl(A_{S_{1,2},S_{4,2}},\{1,3,4\}\bigr)$.
\end{proposition}

\smallskip

Also,
\[
\Pi\subo=\{\beta_1=\epsilon_1-\epsilon_2,\beta_2=\epsilon_2-\epsilon_3,
 \beta_3=\epsilon_3-\epsilon_4,\beta_4=2\epsilon_4\},
\]
with associated Cartan matrix of type $C_4$. The odd part is an
irreducible module of dimension $4\times 8+4\times 2=40$, with
highest weight $R(\epsilon_1+\epsilon_2+\epsilon_3)=\omega_3$.

As for the Lie superalgebra $\frg(S_1,S_{4,2})$, the irreducible
module for the symplectic Lie algebra $\frsp_8$ with this highest
weight is obtained as follows. Let $W$ be the natural eight
dimensional module for $\frsp_8$, and let
$\varphi:\bigwedge^3W\rightarrow W$ be the linear map such that
$\varphi(z_1\wedge z_2\wedge z_3)=\{z_1\vert z_2\}z_3+\{z_2\vert
z_3\}z_1+\{z_3\vert z_1\}z_2$. This time $\ker\varphi$ is not
irreducible but, since the characteristic is $3$, contains the
irreducible submodule $\tilde W=\{\sum_{i=1}^4 a_i\wedge b_i\wedge
z_i :z\in W\}$, which is isomorphic to $W$ (here $\{a_i,b_i:
i=1,\ldots,4\}$ is a symplectic basis of $W$). Then
$\ker\varphi/\tilde W$ is an irreducible module of dimension $40$
and the class of $b_1\wedge b_2\wedge b_3$ modulo $\tilde W$ is a
highest weight vector of weight $\omega_3$.

\begin{proposition}\label{eq:pr:S12S42evenodd}
$\frg(S_{1,2},S_{4,2})$ is a simple Lie superalgebra with even part
isomorphic to the symplectic Lie algebra $\frsp_8$ and odd part
isomorphic to the irreducible module of dimension $40$ above.
\end{proposition}

\bigskip

\subsection{$\frg(S_2,S_{4,2})$}\label{sub:gS2S42}

The symmetric composition superalgebra $S_2$ has a basis
$\{e^+,e^-\}$ with multiplication given by
\[
e^\pm\bullet e^\pm=e^\mp,\qquad e^\pm\bullet e^{\mp}=0,
\]
with norm given by $q(e^\pm)=0$, and $b(e^+,e^-)=1$. The orthogonal
Lie algebra $\frso(S_2,q)$ is spanned by $\phi=\sigma_{e^-,e^+}$,
which satisfies $\phi(e^\pm)=\pm e^\pm$, and its triality Lie
algebra is (see \cite[3.4]{ElduqueMagicII}):
\[
\tri(S_2)=\{(\mu_0\phi,\mu_1\phi,\mu_2\phi):\mu_0,\mu_1,\mu_2\in
k,\, \mu_0+\mu_1+\mu_2=0\}.
\]
Besides,
\[
\begin{split}
t_{e^-,e^+}&=\bigl(\sigma_{e^-,e^+},\frac{1}{2}b(e^-,e^+)1-r_{e^-}l_{e^+},
  \frac{1}{2}b(e^-,e^+)1-l_{e^-}r_{e^+}\bigr)\\
   &=(\phi,\phi,\phi),
\end{split}
\]
(since the characteristic is $3$). Consider the two dimensional
abelian Lie algebra
\[
\frt=\{(\mu_0,\mu_1,\mu_2)\in k^3: \mu_0+\mu_1+\mu_2=0\},
\]
with basis $\{t_1=(1,-1,0),t_2=(0,1,-1)\}$. Then $\tri(S_2)$ is
isomorphic to $\frt$ and the action of $\tri(S_2)$ on each
$\iota_i(S_2)$ becomes:
\[
(\mu_0,\mu_1,\mu_2).\iota_i(e^\pm)=\pm\mu_i\iota_i(e^\pm).
\]
Here,
\[
\frg=\frg(S_2,S_{4,2})=\bigl(\tri(S_2)\oplus\tri(S_{4,2})\bigr)\oplus
 \bigl(\oplus_{i=0}^2\iota_i(S_2\otimes S_{4,2})\bigr),
\]
so, with standard identifications,
\[
\begin{split}
\frg\subo&=\bigl(\frt\oplus\frsp(V_1)\oplus\frsp(V_2)\oplus\frsp(V_3)\bigr)\oplus\\
 &\qquad\qquad\bigl(S_2^0\otimes V_1\otimes V_2\bigr)\oplus
    \bigl(S_2^1\otimes V_2\otimes V_3\bigr)\oplus
    \bigl(S_2^2\otimes V_1\otimes V_3\bigr)\\[6pt]
\frg\subuno&=\bigl(V_1\otimes V_2\otimes V_3\bigr)\oplus
  \bigl(S_2^0\otimes V_3\bigr)\oplus
    \bigl(S_2^1\otimes V_1\bigr)\oplus
    \bigl(S_2^2\otimes V_2\bigr),
\end{split}
\]
where, for any $i=0,1,2$, $S_2^i$ is just a copy of $S_2$, with
basic elements $e_i^{\pm}$.

Consider the free $\bZ$-module $F$ with basis
$\{\delta_1,\delta_2,\epsilon_1,\epsilon_2,\epsilon_3\}$. The Lie
superalgebra $\frg$ is graded over $F\,(\cong \bZ^5)$ by assigning
degree $\epsilon_i$ to $w_i$, $-\epsilon_i$ to $v_i$, $i=1,2,3$, as
in Example \ref{ex:d21}, and
\begin{equation}\label{eq:degreesS2}
\deg(e_0^\pm)=\pm\delta_1,\quad \deg(e_1^\pm)=\pm\delta_2,\quad
 \deg(e_2^\pm)=\mp(\delta_1+\delta_2).
\end{equation}

In this way, the sets of even and odd nonzero degrees that appear in
$\frg(S_2,S_{4,2})$ are:
\[
\begin{split}
\Phi\subo&=\{\pm 2\epsilon_i:i=1,2,3\}\\ &\qquad \cup
 \{\pm\delta_1\pm\epsilon_1\pm\epsilon_2,
 \pm\delta_2\pm\epsilon_2\pm\epsilon_3,
 \pm(\delta_1+\delta_2)\pm\epsilon_1\pm\epsilon_3\},\\[6pt]
\Phi\subuno&=\{\pm\epsilon_1\pm\epsilon_2\pm\epsilon_3,
   \pm\delta_1\pm\epsilon_3,\pm\delta_2\pm\epsilon_1,
   \pm(\delta_1+\delta_2)\pm\epsilon_2\}.
\end{split}
\]
The abelian subalgebra $\frh=kt_1\oplus kt_2\oplus kh_1\oplus
kh_2\oplus kh_3$ ($h_i=\gamma_{v_i,w_i}$, $i=1,2,3$), is a Cartan
subalgebra of $\frg(S_2,S_{4,2})$ and the image of
$\Phi\subo\cup\Phi\subuno$ under the natural homomorphism of abelian
groups
\[
\begin{split}
R:F&\longrightarrow \frh^*\\
\delta_1&\mapsto R(\delta_1)\,(:t_1\mapsto 1,\, t_2\mapsto
  0,\, h_j\mapsto 0),\\
\delta_2&\mapsto R(\delta_2)\,(:t_1\mapsto -1,\, t_2\mapsto
  1,\, h_j\mapsto 0),\\
\epsilon_i&\mapsto R(\epsilon_i)\,(:t_1,t_2\mapsto 0,\,
   h_j\mapsto\delta_{ij}),
\end{split}
\]
is precisely the set of roots of $\frg(S_2,S_{4,2})$ relative to
$\frh$.

Consider the lexicographic order on $F$ with
$\delta_1>\delta_2>\epsilon_1>\epsilon_2>\epsilon_3>0$. Then,
\[
\Phi=\{\alpha_1=\delta_1-\epsilon_1-\epsilon_2,
 \alpha_2=2\epsilon_2,\alpha_3=\epsilon_1-\epsilon_2-\epsilon_3,
 \alpha_4=2\epsilon_3,\alpha_5=\delta_2-\epsilon_1\}
\]
is a linearly independent set with $\Phi\subseteq \bZ\Pi$, so that
any positive element in $\Phi$ is a sum of elements in $\Pi$.
Consider the elements:
\[
\begin{aligned}
E_1&=e_0^+\otimes v_1\otimes v_2,&H_1&=(t_1-t_2)+h_1+h_2,&
  F_1&=\xi_1e_0^-\otimes w_1\otimes w_2,\\
E_2&=\gamma_{w_2,w_2},&H_2&=h_2,&F_2&=\xi_2\gamma_{v_2,v_2},\\
E_3&=w_1\otimes v_2\otimes v_3,&H_3&=-h_1+h_2+h_3,&
   F_3&=\xi_3v_1\otimes w_2\otimes w_3,\\
E_4&=\gamma_{w_3,w_3},&H_4&=h_3,&F_4&=\xi_4\gamma_{v_3,v_3},\\
E_5&=e_1^+\otimes v_1,&H_5&=(t_1-t_2)+h_1,&
  F_5&=\xi_5e_1^-\otimes w_1,
\end{aligned}
\]
where $\xi_1,\ldots,\xi_5$ are suitable scalars so as to have
$[E_i,F_j]=\delta_{ij}H_j$ for any $i,j$.

With these elements $\frg(S_2,S_{4,2})$ is $\bZ$-graded, by assining
degree $1$ to $E_1,\ldots,$ $E_5$, and degree $-1$ to
$F_1,\ldots,F_5$, and the hypotheses of Theorem
\ref{th:centerlessgAtau}  are satisfied relative to the rank $4$
matrix (and Dynkin diagram):
\[
A_{S_2,S_{4,2}}=\begin{pmatrix} 2&-1&0&0&0\\ -1&2&-1&0&0\\
0&-1&0&-1&1\\ 0&0&-1&2&0\\ 0&0&1&0&0\end{pmatrix},\qquad
\SdosScuatrodos
\]

Therefore:

\begin{proposition}\label{pr:gS2S42}
The Lie superalgebra $\frg(S_2,S_{4,2})$ is isomorphic to the
centerless contragredient Lie superalgebra
$\frg\bigl(A_{S_2,S_{4,2}},\{3,5\}\bigr)/\frc$. It is not simple,
but its derived superalgebra $[\frg(S_2,S_{4,2}),\frg(S_2,S_{4,2})]$
is a simple ideal of codimension $1$.
\end{proposition}

\smallskip

Also,
\[
\Pi\subo=\{\beta_1=2\epsilon_1,\beta_2=\delta_1-\epsilon_1-\epsilon_2,
 \beta_3=2\epsilon_2,\beta_4=\delta_2-\epsilon_2-\epsilon_3,
 \beta_5=2\epsilon_3\}
\]
is a linearly independent set with $\Phi\subo\subseteq \bZ\Pi\subo$,
so that any positive element in $\Phi$ is a sum of elements in
$\Pi\subo$, and with similar arguments one obtains that the
associated Cartan matrix is of type $A_5$, $\frg\subo$ thus being
isomorphic to the centerless contragredient Lie algebra
$\frg(A_5)/\frc$, which is isomorphic to the projective general Lie
algebra $\frpgl_6$. The highest degree in $\frg\subuno$ is
$\delta_1+\delta_2+\epsilon_3$, so the highest weight is
$R(\delta_1+\delta_2+\epsilon_3)=\omega _3$. Therefore:

\begin{proposition}\label{pr:gS2S42evenodd}
$\frg(S_2,S_{4,2})$ is a Lie superalgebra with even part isomorphic
to $\frpgl_6$, and odd part isomorphic to the third exterior power
of the natural module for $\frgl_6$ (of dimension $20$).
\end{proposition}

Since the characteristic is $3$, this third exterior power is indeed
a module for $\frpgl_6=\frgl_6/k1$.

\bigskip

\subsection{$\frg(S_{1,2},S_{1,2})$}\label{sub:gS12S12}

Now, let us consider two copies of $S_{1,2}$: $k1\oplus V_i$,
$i=1,2$ (Example \ref{ex:B12bar}). Then the even and odd parts of
\[
\frg=\frg(S_{1,2},S_{1,2})=\bigl(\tri(S_{1,2})\oplus\tri(S_{1,2})\bigr)
 \oplus\bigl(\oplus_{i=0}^2\iota_i(S_{1,2}\otimes S_{1,2})\bigr)
\]
are given, because of Corollary \ref{co:triB12bar}, by:
\[
\begin{split}
\frg\subo&=\bigl(\frsp(V_1)\oplus\frsp(V_2)\bigr)\oplus
 \bigl(\oplus_{i=0}^2k\iota_i(1\otimes 1)\bigr)
 \oplus\bigl(\oplus_{i=0}^2\iota_i(V_1\otimes V_2)\bigr),\\
\frg\subuno&=\Bigl(V_1\oplus
  \bigl(\oplus_{i=0}^2\iota_i(V_1\otimes 1)\bigr)\Bigr)\oplus
    \Bigl(V_2\oplus
  \bigl(\oplus_{i=0}^2\iota_i(1\otimes V_2)\bigr)\Bigr).
\end{split}
\]
Consider the quaternion algebra
\[
Q=k1\oplus\bigl(\oplus_{i=0}^2kx_i\bigr),
\]
with $x_i^2=-1$ for any $i=0,1,2$, and
$x_ix_{i+1}=-x_{i+1}x_i=x_{i+2}$ for any $i$ (indices modulo $3$).
Its norm $N$ satisfies $N(1)=N(x_i)=1$, $i=0,1,2$, and
$N(1,x_i)=N(x_i,x_j)=0$ for any $i\ne j$
($N(a,b)=N(a+b)-N(a)-N(b)$).

Since the characteristic of $k$ is $3$, $N(1+x_1+x_2)=0$, hence $N$
represents $0$ and $Q$ is the split quaternion algebra. That is, $Q$
is isomorphic to $\Mat_2(k)$. In particular,
$Q_0=(k1)^\perp=kx_0\oplus kx_1\oplus kx_2$ is a Lie algebra under
the commutator $[x,y]=xy-yx$, which is isomorphic to $\frsl_2$.

The subspace $\oplus_{i=0}^2k\iota_i(1\otimes 1)$ in $\frg\subo$ is
a subalgebra, and
\[
[\iota_i(1\otimes 1),\iota_{i+1}(1\otimes
1)]=\iota_{i+2}\bigl((1\bullet 1)\otimes (1\bullet 1)\bigr)=
\iota_{i+2}(1\otimes 1),
\]
while $[x_i,x_{i+1}]=2x_{i+2}=-x_{i+2}$. Therefore,
$\oplus_{i=0}^2k\iota_i(1\otimes 1)$ is isomorphic to $Q_0$ under
the linear map that takes $\iota_i(1\otimes 1)$ to $-x_i$.

As vector spaces, there is a natural isomorphism
\[
\Gamma\subo:\frg\subo\rightarrow
 \bigl(\frsp(V_1)\oplus\frsp(V_2)\oplus Q_0\bigr)
 \oplus\bigl(V_1\otimes V_2\otimes Q_0\bigr),
\]
which is the identity on $\frsp(V_1)\oplus \frsp(V_2)$, and such
that
\[
\Gamma\subo\bigl(\iota_i(1\otimes 1)\bigr)=-x_i,\quad
 \Gamma\subo\bigl(\iota_i(u_1\otimes u_2)\bigr)=u_1\otimes
 u_2\otimes x_i,
\]
for any $i=0,1,2$ and $u_1\in V_1$, $u_2\in V_2$.

The Lie bracket in $\frg\subo$ is then transfered to the right hand
side as follows:
\begin{itemize}
\item $\frsp(V_1)\oplus\frsp(V_2)\oplus Q_0$ is a Lie subalgebra
with componentwise bracket.%
\smallskip

\item $[s_1,u_1\otimes u_2\otimes p]=s_1(u_1)\otimes u_2\otimes p$,
$[s_2,u_1\otimes u_2\otimes p]=u_1\otimes s_2(u_2)\otimes p$, for
any $s_i\in \frsp(V_i)$, $u_i\in V_i$, $i=1,2$, and $p\in Q_0$.
\smallskip

\item $[q,u_1\otimes u_2\otimes p]=u_1\otimes u_2\otimes [q,p]$, for
any $u_i\in V_i$, $i=1,2$, and $p,q\in Q_0$. This is because
$[\iota_i(1\otimes 1),\iota_i(u_1\otimes u_2)]=0$, and
$[\iota_i(1\otimes 1),\iota_{i+1}(u_1\otimes
u_2)]=\iota_{i+2}(u_1\otimes u_2)$ and $[\iota_i(1\otimes
1),\iota_{i+2}(u_1\otimes u_2)]=-\iota_{i+1}(u_1\otimes u_2)$ in
$\frg\subo$ (see Section 3).
\smallskip

\item For any $u_i,u_i'\in V_i$, $i=1,2$, and $p,p'\in Q_0$:
\begin{multline*}
[u_1\otimes u_2\otimes p,u_1'\otimes u_2'\otimes p']\\ =
  -\langle u_2\vert u_2'\rangle N(p,p')\gamma_{u_1,u_1'}
   -\langle u_1\vert u_1'\rangle N(p,p')\gamma_{u_2,u_2'}-
   \langle u_1\vert u_1'\rangle\langle u_2\vert u_2'\rangle [p,q].
\end{multline*}
This is because $[\iota_i(u_1\otimes u_2),\iota_i(u_1'\otimes
u_2')]=\langle u_2\vert u_2'\rangle \gamma_{u_1,u_1'}+\langle
u_1\vert u_1'\rangle \gamma_{u_2,u_2'}$ (note that $S_{1,2}=k1\oplus
V$, and $\sigma_{u,u'}=-\gamma_{u,u'}$ \eqref{eq:sigmauvB12}),
$N(1,1)=2=-1$ and $[\iota_i(u_1\otimes u_2),\iota_{i+1}(u_1'\otimes
u_2')]=-\iota_{i+2}\bigl((u_1\bullet u_1')\otimes (u_2\bullet
u_2')\bigr)=-\langle u_1\vert u_1'\rangle\langle u_2\vert
u_2'\rangle\iota_{i+2}(1\otimes 1)$, which corresponds to
$[u_1\otimes u_2\otimes x_i,u_1'\otimes u_2'\otimes x_{i+1}]=
\langle u_1\vert u_1'\rangle\langle u_2\vert u_2'\rangle x_{i+2}$
($[x_i,x_{i+1}]=-x_{i+2}$).
\end{itemize}

\medskip

Now, as vector spaces, there is also a natural isomorphism:
\[
\begin{split}
\Gamma\subuno: \frg\subuno&\longrightarrow
  (V_1\oplus V_2)\otimes Q,\\
 u_1&\mapsto\ u_1\otimes 1,\\
 u_2&\mapsto\ -u_2\otimes 1,\\
  \iota_i(u_1\otimes 1)&\mapsto -u_1\otimes x_i,\\
  \iota_i(1\otimes u_2)&\mapsto u_2\otimes x_i,
\end{split}
\]
for $i=0,1,2$, $u_1\in V_1$, and $u_2\in V_2$. Under $\Gamma\subo$
and $\Gamma\subuno$, the Lie bracket $\frg\subo\times
\frg\subuno\rightarrow \frg\subuno$ is transferred to:

\begin{itemize}
\item $[s_i,(u_1+u_2)\otimes q]=s_i(u_i)\otimes q$ for any $s_i\in
\frsp(V_i)$, $u_i\in V_i$, $i=1,2$ and $q\in Q$.
\smallskip

\item $[p,(u_1+u_2)\otimes q]=(u_1+u_2)\otimes pq$, for any $p\in
Q_0$, $q\in Q$, and $u_i\in V_i$, $i=1,2$.
\smallskip

\item $[u_1\otimes u_2\otimes p,u_1'\otimes q]=-\langle u_1\vert
u_1'\rangle u_2\otimes pq$, for any $u_1,u_1'\in V_1$, $u_2\in V_2$,
$p\in Q_0$ and $q\in Q$.
\smallskip

\item $[u_1\otimes u_2\otimes p,u_2'\otimes q]=\langle u_2\vert
u_2'\rangle u_1\otimes pq$, for any $u_1\in V_1$, $u_2,u_2'\in V_2$,
$p\in Q_0$ and $q\in Q$.

\end{itemize}
\medskip

While the bracket $\frg\subuno\times \frg\subuno\rightarrow
\frg\subo$ transfers to:

\begin{itemize}
\item $[u_1\otimes q,u_1'\otimes q']=-N(q,q')\gamma_{u_1,u_1'}
 -\langle u_1,u_1'\rangle (q\bar q'-q'\bar q)$,
\smallskip

\item $[u_2\otimes q,u_2'\otimes q']=-N(q,q')\gamma_{u_2,u_2'}
 -\langle u_2,u_2'\rangle (q\bar q'-q'\bar q)$,
\smallskip

\item $[u_1\otimes q,u_2\otimes q']=-u_1\otimes u_2\otimes (q\bar
q'-q'\bar q)$,
\end{itemize}
\medskip

\noindent for any $u_i,u_i'\in V_i$, $i=1,2$, and $q,q'\in Q$. Here
$\bar{\phantom{q}}$ denotes the canonical involution in the
quaternion algebra $Q$ (that is, the symplectic involution in
$\Mat_2(k)$).

Therefore, the even and odd parts of $\frg=\frg(S_{1,2},S_{1,2})$
can be identified to
\begin{equation}\label{eq:gS12S12g0g1}
\begin{split}
\frg\subo&=\bigl(\frsp(V_1)\oplus\frsp(V_2)\oplus\frsl_2\bigr)
   \oplus\bigl(V_1\otimes V_2\otimes \frsl_2\bigr),\\
\frg\subuno&=(V_1\oplus V_2)\otimes \frgl_2.
\end{split}
\end{equation}

The action of $\frsl_2$ on $\frgl_2$ is given by left
multiplication, hence $\frgl_2$ decomposes, as a module for
$\frsl_2$, as $\frgl_2=\frj_1\oplus \frj_{-1}$, where $\frj_1$
(respectively $\frj_{-1}$) consists of the $2\times 2$ matrices with
zero second (resp. first) column.

Let $\{v_i,w_i\}$ be a symplectic basis of $V_i$, $i=1,2$, as
considered so far, and let $\{h=E_{11}-E_{22},e=E_{12},f=E_{21}\}$
be the standard basis of $\frsl_2$ ($E_{ij}$ denotes the $2\times 2$
matrix with entry $(i,j)$ equal to $1$ and all the other entries
$0$). Note that $\{E_{11},E_{21}\}$ is a basis of $\frj_1$, while
$\{ E_{12},E_{22}\}$ is a basis of $\frj_{-1}$.

Consider here the free $\bZ$-module $F$ with basis
$\epsilon_1,\epsilon_2,\epsilon,\delta$. The Lie superalgebra
$\frg=\frg(S_{1,2},S_{1,2})$ is $F$-graded by assigning  degree
$\epsilon_i$ to $w_i$ and $-\epsilon_i$ to $v_i$ as usual, $i=1,2$,
degree $\epsilon+\delta$ to $E_{11}\,(\in\frj_1)$,
$-\epsilon+\delta$ to $E_{21}\,(\in\frj_1)$, $\epsilon-\delta$ to
$E_{12}\,(\in\frj_{-1})$, and $-\epsilon-\delta$ to
$E_{22}\,(\in\frj_{-1})$. Then it follows that, for instance, the
degree of $v_1\otimes E_{22}$ is $-\epsilon_1-\epsilon-\delta$, or
that the degree of $e=E_{12}\in\frsl_2$ is $2\epsilon$.

The sets of nonzero even and odd degrees are:
\[
\begin{split}
\Phi\subo&=\{\pm 2\epsilon_1,\pm 2\epsilon_2,\pm 2\epsilon,
  \pm\epsilon_1\pm\epsilon_2,\pm\epsilon_1\pm\epsilon_2\pm
  2\epsilon\},\\
\Phi\subuno&=\{ \pm\epsilon_1\pm\epsilon\pm\delta,
   \pm\epsilon_2\pm\epsilon\pm\delta\}.
\end{split}
\]

The abelian subalgebra $\frh=kh_1\oplus kh_2\oplus kh$
($h_i=\gamma_{v_i,w_i}\in \frsp(V_i)$, $i=1,2$, as usual) is a
Cartan subalgebra of $\frg$ and the image of
$\Phi\subo\cup\Phi\subuno$ under the homomorphism of abelian groups:
\[
\begin{split}
R:F&\longrightarrow \frh^*\\
 \epsilon_1&\mapsto R(\epsilon_1)\,(:h_1\mapsto 1,\, h_2\mapsto 0,\,
 h\mapsto 0),\\
 \epsilon_2&\mapsto R(\epsilon_2)\,(:h_1\mapsto 0,\, h_2\mapsto 1,\,
 h\mapsto 0),\\
 \epsilon&\mapsto R(\epsilon)\,(:h_1\mapsto 0,\, h_2\mapsto 0,\,
 h\mapsto 1),\\
 \delta&\mapsto 0,
\end{split}
\]
is precisely the set of roots of $\frg$ relative to $\frh$.

\smallskip

Consider the lexicographic order on $F$ with
$\delta>\epsilon_1>\epsilon_2>\epsilon>0$, then
\[
\Phi=\{\alpha_1=2\epsilon_2,\alpha_2=\epsilon_1-\epsilon_2-2\epsilon,
 \alpha_3=2\epsilon,\alpha_4=\delta-\epsilon_1-\epsilon\}
\]
is a linearly independent set in $F$ with
$\Phi=\Phi\subo\cup\Phi\subuno\subseteq \bZ\Pi$, so that any
positive element in $\Phi$ is a sum of elements in $\Pi$. Consider
the elements:
\[
\begin{aligned}
E_1&=\gamma_{w_2,w_2},&E_2&=w_1\otimes v_2\otimes E_{21},&
 E_3&=e,&E_4&=v_1\otimes E_{21},\\
F_1&=-\gamma_{v_2,v_2},&F_2&=-v_1\otimes w_2\otimes E_{12},&
 F_3&=f,&F_4&=w_1\otimes E_{12}.
\end{aligned}
\]
Then:
\[
\begin{aligned}
H_1=[E_1,F_1]&=\gamma_{v_2,w_2}=h_2,\\
 H_2=[E_2,F_2]&=-\gamma_{v_1,w_1}+\gamma_{v_2,w_2}+h=-h_1+h_2+h,\\
H_3=[E_3,F_3]&=h,\\
 H_4=[E_4,F_4]&=\gamma_{v_1,w_1}-h=h_1-h,
\end{aligned}
\]
and the hypotheses of Theorem \ref{th:centerlessderivedgAtau} are
satisfied relative to the rank $3$ matrix (and Dynkin diagram):
\[
A_{S_{1,2},S_{1,2}}=\begin{pmatrix}
 2&-1&0&0\\ -1&2&-1&0\\ 0&-2&2&-1\\ 0&0&1&0\end{pmatrix},\qquad
 \SunodosSunodos
\]

Therefore:

\begin{proposition}\label{pr:gS12S12}
The Lie superalgebra $\frg(S_{1,2},S_{1,2})$ is isomorphic to the
centerless derived contragredient Lie superalgebra
$\frg'\bigl(A_{S_{1,2},S_{1,2}},\{4\}\bigr)/\frc$.
\end{proposition}

\smallskip

Moreover, $\{E_i,F_i,H_i: i=1,2,3\}$ generate the even part of
$\frg\subo$, the associated Cartan subalgebra is of type $B_3$,
$\frg\subo$ thus being isomorphic to the orthogonal Lie algebra
$\frso_7$. The grading of $\frg$ by $F$ gives, in particular, a
consistent $\bZ$-grading by means of the projection $F=\bZ
\epsilon_1\oplus\bZ\epsilon_2\oplus\bZ\epsilon\oplus\bZ\delta\rightarrow
\bZ$, $n_1\epsilon_1+n_2\epsilon_2+n_3\epsilon+n_4\delta\mapsto
n_4$, where
\[
\frg=\frg_{-1}\oplus\frg_0\oplus\frg_1,
\]
with
\[
\frg_0=\frg\subo,\quad \frg_1=(V_1\oplus V_2)\otimes \frj_1,\quad
 \frg_{-1}=(V_1\oplus V_2)\otimes \frj_{-1}.
\]
As modules over $\frg\subo=\frg_0$, $\frg_1$ and $\frg_{-1}$ are
isomorphic, and $\frg_1$ is an irreducible module with highest
weight $R(\epsilon_1+\epsilon+\delta)=\omega_3$, so $\frg_1$ is the
spin module for $\frg_0$. Hence:

\begin{proposition}\label{pr:gS12S12evenodd}
$\frg(S_{1,2},S_{1,2})$ is a simple Lie superalgebra whose even part
is isomorphic to $\frso_7$, and its odd part is the direct sum of
two copies of the spin module for the even part.
\end{proposition}

\begin{corollary}\label{co:S12S12}
The Lie superalgebra $\frg(S_{1,2},S_{1,2})$ is isomorphic to the
Lie superalgebra in \cite[Theorem 4.23(ii)]{ElduqueNewSimple3},
constructed in terms of a null orthogonal triple system.
\end{corollary}

\bigskip

\subsection{$\frg(S_1,S_{1,2})$}

The Lie superalgebra $\frg(S_1,S_{1,2})$ is a subsuperalgebra of
$\frg(S_{1,2},S_{1,2})$. Given the description of this latter
superalgebra in \ref{sub:gS12S12}, more specifically in
\eqref{eq:gS12S12g0g1}, the even and odd parts of
$\frg=\frg(S_1,S_{1,2})$ can be described as
\[
\begin{split}
\frg\subo&=\frsp(V_2)\oplus\frsl_2\ (\subseteq
    \frg(S_{1,2},S_{1,2})\subo),\\
\frg\subuno&= V_2\otimes \frgl_2\ (\subseteq
           \frg(S_{1,2},S_{1,2})\subuno),
\end{split}
\]
with
\[
\Phi\subo=\{\pm 2\epsilon_2,\pm 2\epsilon\},\quad
 \Phi\subuno=\{ \pm\epsilon_2\pm\epsilon\pm\delta\}.
\]
With the lexicographic order in which
$\delta>\epsilon_2>\epsilon>0$,
\[
\Pi=\{\alpha_1=2\epsilon_2,\alpha_2=\delta-\epsilon_2-\epsilon,
\alpha_3=2\epsilon\},
\]
with associated matrix and Dynkin diagram:
\[
A_{S_1,S_{1,2}}=\begin{pmatrix}
  2&-1&0\\ -1&0&1\\ 0&-1&2 \end{pmatrix},\qquad \SunoSunodos
\]
and the same arguments as in the previous section give:

\begin{proposition}\label{pr:gS1S12}
The Lie superalgebra $\frg(S_1,S_{1,2})$ is isomorphic to the
centerless derived contragredient Lie superalgebra
$\frg'\bigl(A_{S_1,S_{1,2}},\{2\}\bigr)/\frc$.
\end{proposition}

\begin{corollary}
$\frg(S_1,S_{1,2})$ is isomorphic to the projective special Lie
superalgebra $\frpsl(2,2)$.
\end{corollary}

This is the only Lie superalgebra in the Extended Freudenthal's
Magic Square with a counterpart in characteristic $0$.

Note that the even part $\frg\subo$ is the direct sum of two copies
of $\frsl_2$, while the odd part is a direct sum of two copies of
the four dimensional irreducible module for $\frg\subo$ which
consists of the tensor product of the natural modules for each copy
of $\frsl_2$.

\bigskip

\subsection{$\frg(S_4,S_{1,2})$}

Here a description in terms of $V(\sigma)$'s is again possible.
Three copies of $V$ are needed for $S_4$ and an extra copy for
$S_{1,2}$. The indices $1$, $2$ and $3$ will be used for the three
copies related to $S_4$ and the index $4$ for the extra copy
attached to $S_{1,2}$. Then
\[
\begin{split}
\frg(S_4,S_{1,2})&=\bigl(\tri(S_4)\oplus\tri(S_{1,2})\bigr)\oplus
    \bigl(\oplus_{i=0}^2\iota_i(S_4\otimes S_{1,2})\bigr)\\
    &=\oplus_{\sigma\in\calS_{S_4,S_{1,2}}}V(\sigma),
\end{split}
\]
where
\[
\calS_{S_4,S_{1,2}}=\bigl\{\emptyset,\{4\},\{1,2\},\{1,2,4\},\{2,3\},
 \{2,3,4\},\{1,3\},\{1,3,4\}\bigr\}.
\]
The nonzero even and odd degrees are:
\[
\begin{split}
\Phi\subo&=\{\pm 2\epsilon_i: 1\leq i\leq 4\}\cup
  \{\pm\epsilon_1\pm\epsilon_2,\pm\epsilon_2\pm\epsilon_3,
  \pm\epsilon_1\pm\epsilon_3\},\\
\Phi\subuno&=\{\pm\epsilon_4,\pm\epsilon_1\pm\epsilon_2\pm\epsilon_4,
          \pm\epsilon_2\pm\epsilon_3\pm\epsilon_4,
           \pm\epsilon_1\pm\epsilon_3\pm\epsilon_4\}.
\end{split}
\]
With the lexicographic order given by
$0<\epsilon_1<\epsilon_2<\epsilon_3<\epsilon_4$, the set of
irreducible degrees is
\[
\Pi=\{\alpha_1=\epsilon_3-\epsilon_2,\alpha_2=\epsilon_2-\epsilon_1,
 \alpha_3=2\epsilon_1,\alpha_4=\epsilon_4-\epsilon_2-\epsilon_3\},
\]
which is a linearly independent set over $\bZ$, and
$\Phi=\Phi\subo\cup\Phi\subuno$ is contained in $\bZ\Pi$. The
associated matrix and Dynkin diagram are:
\[
A_{S_4,S_{1,2}}=\begin{pmatrix} 2&-1&0&0\\ -1&2&-2&-1\\
   0&-1&2&0\\ 0&1&0&0\end{pmatrix},\qquad \ScuatroSunodos
\]
and

\begin{proposition}\label{pr:gS4S12}
The Lie superalgebra $\frg(S_4,S_{1,2})$ is isomorphic to the
contragredient Lie superalgebra
$\frg\bigl(A_{S_4,S_{1,2}},\{4\}\bigr)$.
\end{proposition}

\smallskip

Also, the set of irreducible even degrees is
\[
 \Pi\subo=\{
  \alpha_1,\alpha_2,\alpha_3,\tilde\alpha_4=2\epsilon_4\},
\]
with associated Cartan matrix of type $C_3\oplus A_1$. (This also
follows from $\frg(S_4,S_{1,2})\subo=\frg(S_4,S_1)\oplus
\frsp(V_4)$.)

The odd part, which has dimension $26$, is irreducible with highest
weight
$R(\epsilon_2+\epsilon_3+\epsilon_4)=R(\epsilon_2+\epsilon_3)+R(\epsilon_4)$,
so that:

\begin{proposition}\label{pr:gS4S12evenodd}
$\frg(S_4,S_{1,2})$ is a simple Lie superalgebra with even part
isomorphic to the direct sum $\frsp_6\oplus\frsl_2$, and odd part
the irreducible module which is the tensor product of the
irreducible module $V(\omega_2)$ for $\frsp_6$ (of dimension $13$) ,
and the natural two dimensional module for $\frsl_2$.
\end{proposition}

\smallskip

The irreducible module $V(\omega_2)$ for $\frsp_6$ can be described
as follows. Let $W$ be the natural six-dimensional module for
$\frsp_6$, which is endowed with a nondegenerate alternating form
$\{.\vert.\}$, and take a symplectic basis $\{a_i,b_i: i=1,2,3\}$.
Consider the linear map $\varphi:\bigwedge^2W,\rightarrow k$,
$z_1\wedge z_2\mapsto \{z_1\vert z_2\}$. Then, since the
characteristic of $k$ is $3$, $k(\sum_{i=1}^3a_i\wedge b_i)$ is a
trivial submodule of $\ker\varphi$, and
$\ker\varphi/k(\sum_{i=1}^3a_i\wedge b_i)$ is an irreducible module
with highest weight vector the class of $b_2\wedge b_3$ modulo
$k(\sum_{i=1}^3a_i\wedge b_i)$ (conventions as in Subsection
\ref{sub:gS1S42}). This is, up to isomorphism, the required
irreducible module.

\bigskip

\subsection{$\frg(S_8,S_{1,2})$}

Here the indices $1$, $2$, $3$ and $4$ will be used for the copies
of $V$ attached to $S_8$, while the index $5$ will be reserved for
the copy of $V$ attached to $S_{1,2}$. With the same arguments as
before:
\[
\begin{split}
\frg(S_8,S_{1,2})&=\bigl(\tri(S_8)\oplus\tri(S_{1,2})\bigr)\oplus
    \bigl(\oplus_{i=0}^2\iota_i(S_8\otimes S_{1,2})\bigr)\\
    &=\oplus_{\sigma\in\calS_{S_8,S_{1,2}}}V(\sigma),
\end{split}
\]
where
\[
\begin{split}
\calS_{S_8,S_{1,2}}=\bigl\{\emptyset,&\{1,2,3,4\},\{5\},\\
       &\{1,2\},\{3,4\},\{1,2,5\},\{3,4,5\},\\
       &\{2,3\},\{1,4\},\{2,3,5\},\{1,4,5\},\\
       &\{1,3\},\{2,4\},\{1,3,5\},\{2,4,5\}\bigr\}.
\end{split}
\]
Here,
\[
\begin{split}
\Phi\subo&=\{\pm 2\epsilon_i: 1\leq i\leq 5\}\cup
  \{\pm\epsilon_1\pm\epsilon_2\pm\epsilon_3\pm\epsilon_4\}\cup
    \{\pm\epsilon_i\pm\epsilon_j:1\leq i<j\leq 4\},\\
\Phi\subuno&=\{\pm\epsilon_5\}\cup\{\pm\epsilon_i\pm\epsilon_j\pm\epsilon_5:
          1\leq i<j\leq 4\},
\end{split}
\]
and with the lexicographic order given by
$0<\epsilon_1<\epsilon_2<\epsilon_3<\epsilon_4<\epsilon_5$, the set
of irreducible degrees is
\begin{multline*}
\Pi=\{\alpha_1=\epsilon_4-\epsilon_3-\epsilon_2-\epsilon_1,\alpha_2=\epsilon_1,
 \alpha_3=\epsilon_2-\epsilon_1,\\
 \alpha_4=\epsilon_3-\epsilon_2,\alpha_5=\epsilon_5-\epsilon_3-\epsilon_4\},
\end{multline*}
which is a $\bZ$-linearly independent set with $\Phi\subseteq
\bZ\Pi$. The associated matrix and Dynkin diagrams are
\[
A_{S_8,S_{1,2}}=\begin{pmatrix}
 2&-1&0&0&0\\ -1&2&-1&0&0\\ 0&-2&2&-1&0\\
  0&0&-1&2&-1\\ 0&0&0&1&0\end{pmatrix},\qquad \SochoSunodos
\]
and then:

\begin{proposition}\label{pr:gS8S12}
The Lie superalgebra $\frg(S_8,S_{1,2})$ is isomorphic to the
contragredient Lie superalgebra
$\frg\bigl(A_{S_8,S_{1,2}},\{5\}\bigr)$.
\end{proposition}

\smallskip

Also,
\[
\Pi\subo=\{\alpha_1,\alpha_2,\alpha_3,\alpha_4,\tilde\alpha_5=2\epsilon_5\},
\]
with associated Cartan matrix of type $F_4\oplus A_1$ (this also
follows from $\frg(S_8,S_{1,2})\subo=\frg(S_8,S_1)\oplus
\frsp(V_5)$). The odd part, which has dimension $50$, is irreducible
with highest weight
$R(\epsilon_5+\epsilon_4+\epsilon_3)=R(\epsilon_3+\epsilon_4)+R(\epsilon_5)$.
$R(\epsilon_5)$ is the highest weight of the natural two dimensional
module for $\frsp(V_5)\cong\frsl_2$, while
$R(\epsilon_3+\epsilon_4)$ is the highest weight $\omega_4$ for the
direct summand of $\frg(S_8,S_{1,2})\subo$ isomorphic to $\frf_4$,
the exceptional simple split Lie algebra of type $F_4$ (with the
ordering of the roots given above). Thus:

\begin{proposition}\label{pr:gS8S12evenodd}
$\frg(S_8,S_{1,2})$ is a simple Lie superalgebra with even part
isomorphic to $\frf_4\oplus \frsl_2$, and odd part the irreducible
module which is the tensor product of the natural two dimensional
module for $\frsl_2$ and the $25$-dimensional module of highest
weight $\omega_4$ for $\frf_4$.
\end{proposition}

In characteristic $>3$, the module $V(\omega_4)$ for $\frf_4$ is
usually described as the set of trace zero elements $J_0$ in the
split exceptional Jordan algebra $J=H_3(C)$, where $C$ is the split
octonion algebra. However, in characteristic $3$, the identity
matrix belongs to $J_0$, and $J_0/k1$ is the irreducible module
$V(\omega_4)$ in this case.

\bigskip

\subsection{$\frg(S_2,S_{1,2})$}

Take the basis $\{e^+,e^-\}$ of the symmetric composition algebra
$S_2$ as in \ref{sub:gS2S42}, and identify $\tri(S_2)$ with
$\frt=\espan{t_1,t_2}$ as was done there. Take also three copies
$S_2^i$ (with basis $\{e_i^+,e_i^-\}$) of $S_2$, $i=0,1,2$. Here we
consider the Lie superalgebra
\[
\frg=\frg(S_2,S_{1,2})=\bigl(\tri(S_2)\oplus\tri(S_{1,2})\bigr)
  \oplus\bigl(\oplus_{i=0}^2\iota_i(S_2\otimes S_{1,2})\bigr),
\]
so, with standard identifications:
\[
\begin{split}
\frg\subo&=\bigl(\frt\oplus\frsp(V)\bigr)
     \oplus\bigl(\oplus_{i=0}^2 S_2^i\bigr),\\
\frg\subuno&=V\oplus\bigl(\oplus_{i=0}^2(S_2^i\otimes V)\bigr).
\end{split}
\]
Consider the free $\bZ$-module $F$ freely generated by $\delta_1$,
$\delta_2$ and $\epsilon$. With the conventions in \ref{sub:gS2S42},
$\frg$ is $F$-graded and the sets of nonzero even and odd degrees
are:
\[
\begin{split}
\Phi\subo&=\{\pm 2\epsilon,\pm\delta_1,
  \pm\delta_2,\pm(\delta_1+\delta_2)\}\\
\Phi\subuno&=\{\pm\epsilon,\pm\delta_1\pm\epsilon,
 \pm\delta_2\pm\epsilon,\pm(\delta_1+\delta_2)\pm\epsilon\}.
\end{split}
\]
With the lexicographic order in which
$\epsilon>\delta_1>\delta_2>0$,
\[
\Pi=\{\alpha_1=\delta_1,\alpha_2=\delta_2,\alpha_3=\epsilon-(\delta_1+\delta_2)\}
\]
is a linearly independent set with
$\Phi=\Phi\subo\cup\Phi\subuno\subseteq \bZ\Pi$. Now consider the
elements ($\{v,w\}$ denotes a symplectic basis of $V$):
\[
\begin{aligned}
E_1&=e_0^+(=\iota_0(1\otimes e^+)),&E_2&=e_1^+,
 &E_3&=e_2^+\otimes w(=\iota_2(e^+\otimes w)),\\
H_1&=2(t_1-t_2),&H_2&=2(t_1-t_2),&H_3&=(t_1-t_2)-\gamma_{v,w},\\
F_1&=\xi_1e_0^-,&F_2&=\xi_2e_1^-,
 &F_3&=\xi_3e_2^-\otimes v,
\end{aligned}
\]
where $\xi_1,\xi_2,\xi_3$ are suitable scalars so as to have
$[E_i,F_j]=\delta_{ij}H_i$ for any $i,j$.

With these elements, $\frg(S_2,S_{1,2})$ is $\bZ$-graded, by
assigning degree $1$ to $E_1,E_2,E_3$ and degree $-1$ to
$F_1,F_2,F_3$, and the hypotheses of Theorem
\ref{th:centerlessgAtau} are satisfied relative to the rank $2$
matrix (with associated Dynkin diagram):
\[
A_{S_2,S_{1,2}}=\begin{pmatrix}
 2&-1&-1\\ -1&2&-1\\ 1&1&0\end{pmatrix},\qquad\qquad \SdosSunodos
\]

Therefore:

\begin{proposition}\label{pr:gS2S12}
The Lie superalgebra $\frg(S_2,S_{1,2})$ is isomorphic to the
centerless contragredient Lie superalgebra
$\frg\bigl(A_{S_2,S_{1,2}},\{3\}\bigr)/\frc$. It is not simple, but
its derived superalgebra is a simple ideal of codimension $1$.
\end{proposition}

\smallskip

Also,
\[
\Pi\subo=\{\beta_1=\delta_1,\beta_2=\delta_2,\beta_3=2\epsilon\}
\]
is a linearly independent set with $\Phi\subo\subseteq \bZ\Pi\subo$
and, with already used arguments, it is obtained that the associated
Cartan matrix is of type $A_2\oplus A_1$, $\frg\subo$ thus being
isomorphic to the direct sum of the centerless contragredient Lie
algebra $\frg(A_2)/\frc$, which is isomorphic to $\frpgl_3$, and of
$\frsl_2$. The highest degree in $\frg\subuno$ is
$(\delta_1+\delta_2)+\epsilon$, and $\frg\subuno$ is easily seen to
be irreducible. Hence:

\begin{proposition}\label{pr:gS2S12evenodd}
$\frg(S_2,S_{1,2})$ is a Lie superalgebra with even part isomorphic
to the direct sum of $\frpgl_3$ and $\frsl_2$, and odd part
isomorphic, as a module for the even part, to the tensor product of
the ``adjoint module'' $\frpsl_3$ of $\frpgl_3$ by the natural two
dimensional module for $\frsl_2$.
\end{proposition}

\bigskip

Table \ref{ta:summary} summarizes the information on the even and
odd parts of the Lie superalgebras in the extended Freudenthal Magic
Square obtained in this section. The notation $(n)$ will indicate a
module of dimension $n$.

\begin{table}[h!]
$$
\vbox{\offinterlineskip
 \halign{\hfil$#$\quad\hfil&%
 \vrule height 12pt width1pt depth 6pt #%
 &\hfil\quad$#$\quad\hfil%
 &\vrule  depth 6pt width .5pt #%
 &\hfil\qquad$#$\hfil\quad\cr
 \bigstrut &width 0pt&S_{1,2}&\omit\vrule height 8pt depth 6pt width .5pt&S_{4,2}\cr
 &\multispan4{\hreglonfill}\cr
 S_1&&\frpsl_{2,2}&&\frsp_6\oplus (14)\cr
 \bigstrut S_2&&
   \bigl(\frpgl_3\oplus\frsl_2\bigr)\oplus\bigl(\frpsl_3\otimes (2)\bigr)&&
    \frpgl_6\oplus (20)\cr
 \bigstrut S_4&&
   \bigl(\frsp_6\oplus\frsl_2\bigr)\oplus\bigl((13)\otimes (2)\bigr)
    &&\frso_{12}\oplus spin_{12}\cr
 \bigstrut S_8&&
   \bigl(\frf_4\oplus\frsl_2\bigr)\oplus\bigl((25)\otimes (2)\bigr)&&
      \fre_7\oplus (56)\cr
 \multispan5{\hregletafill}\cr
 \bigstrut S_{1,2}&&\frso_7\oplus 2spin_7 &&\frsp_8\oplus(40)\cr
 \bigstrut S_{4,2}&& & & \frso_{13}\oplus spin_{13}\cr}}
$$
\bigskip
\caption{Even and odd parts}\label{ta:summary}
\end{table}

\medskip

With the exception of $\frpsl_{2,2}$, none of these Lie
superalgebras have a counterpart in characteristic $0$. Corollary
\ref{co:S42S42} shows that $\frg(S_{4,2},S_{4,2})$ has already
appeared in \cite[Theorem 3.1(ii)]{ElduqueNewSimple35}, while
Corollary \ref{co:S12S12} shows that $\frg(S_{1,2},S_{1,2})$ is the
Lie superalgebra that appears in \cite[Theorem
4.23(ii)]{ElduqueNewSimple3}.

In a forthcoming paper, it will be shown that the Lie superalgebras
$\frg(S_r,S_{1,2})$ and $\frg(S_r,S_{4,2})$, for $r=1$, $4$ and $8$,
and their derived subalgebras for $r=2$, are precisely the simple
Lie superalgebras defined in \cite{ElduqueNewSimple3} in terms of
symplectic and orthogonal triple systems related to simple Jordan
algebras of degree $3$.



\providecommand{\bysame}{\leavevmode\hbox
to3em{\hrulefill}\thinspace}
\providecommand{\MR}{\relax\ifhmode\unskip\space\fi MR }
\providecommand{\MRhref}[2]{%
  \href{http://www.ams.org/mathscinet-getitem?mr=#1}{#2}
} \providecommand{\href}[2]{#2}

\end{document}